\documentclass[12pt, reqno]{amsart}
\usepackage{pifont}
\usepackage{mathrsfs}
\usepackage{geometry}
\usepackage{titletoc}
\usepackage{stix2}

\usepackage{amsmath}
\usepackage{amssymb} 
\usepackage{enumitem} 
\usepackage{mathtools} 
\usepackage[table]{xcolor} 
\usepackage[all]{xy} 
\usepackage{tikz} 
\usepackage{tikz-cd}
\usepackage{indentfirst} 
\usepackage{babel} 
\usepackage{setspace} 

\usepackage[colorlinks,linkcolor=red,anchorcolor=green,citecolor=blue]{hyperref} 
\hypersetup{linktocpage = true} 

\usepackage{rotating} 

\usepackage{ytableau} 
\usepackage{longtable} 
\newcolumntype{M}[1]{>{\centering\arraybackslash}m{#1}} 

\usepackage{fancyhdr}

\DeclareMathOperator*{\Sym}{Sym}

\newcommand\bp{{\bar\partial}}

\theoremstyle{plain}
\newtheorem{thm}{Theorem}[section]
\newtheorem{lemma}[thm]{Lemma}
\newtheorem{prop}[thm]{Proposition}
\newtheorem{cor}[thm]{Corollary}
\newtheorem{defn}[thm]{Definition}

\newtheorem{conjecture}[thm]{Conjecture}

\theoremstyle{definition}
\newtheorem{example}[thm]{Example}
\newtheorem{remark}[thm]{Remark}

\newcommand{\btheorem}{\begin{thm}}
    \newcommand{\etheorem}{\end{thm}}
\newcommand{\bproposition}{\begin{prop}}
    \newcommand{\eproposition}{\end{prop}}
\newcommand{\bdefinition}{\begin{defn}}
    \newcommand{\edefinition}{\end{defn}}
\newcommand{\bcorollary}{\begin{cor}}
    \newcommand{\ecorollary}{\end{cor}}
\newcommand{\bproof}{\begin{proof}}
    \newcommand{\eproof}{\end{proof}}
\newcommand{\bremark}{\begin{remark}}
    \newcommand{\eremark}{\end{remark}}
\newcommand{\eexample}{\end{example}}
\newcommand{\bexample}{\begin{example}}

\newcommand{\elemma}{\end{lemma}}
\newcommand{\blemma}{\begin{lemma}}

\newcommand{\la}{\langle}
\newcommand{\ra}{\rangle}
\newcommand{\sq}{\sqrt{-1}}

\newcommand{\p}{\partial}

\renewcommand{\bar}{\overline}

\renewcommand{\phi}{\varphi}

\newcommand{\beq}{\begin{equation}}
\newcommand{\eeq}{\end{equation}}
\newcommand{\ee}{\end{eqnarray*}}
\newcommand{\be}{\begin{eqnarray*}}

\newcommand{\bd}{\begin{enumerate}}
    \newcommand{\ed}{\end{enumerate}}

\renewcommand{\hat}{\widehat}
\renewcommand{\tilde}{\widetilde}

\newcommand{\qtq}[1]{\quad\mbox{#1}\quad}
\renewcommand{\bp}{\bar{\partial}}
\newcommand{\Om}{\Omega}

\newcommand{\ts}{\otimes}

\renewcommand{\>}{\rightarrow}

\newcommand{\C}{{\mathbb C}}

\renewcommand{\P}{{\mathbb P}}

\newcommand{\R}{{\mathbb R}}

\usepackage{harmony}

\newcommand{\vone}{\vskip 1\baselineskip}

\renewcommand{\#}{\sharp}

\newcommand{\nm}[1]{\left\Vert #1\right\Vert}
\setlist[itemize]{leftmargin=*}
\setlist[enumerate]{leftmargin=*}

\numberwithin{equation}{section} 

\makeatletter

\usepackage{fancyhdr}
\renewcommand{\leftmark}{{\scriptsize Weitzenb\"ock-Bochner-Kodaira formulas with quadratic curvature terms}}

\title{Weitzenb\"ock-Bochner-Kodaira formulas with quadratic curvature terms}

\author{Mingwei Wang}
\author{Xiaokui Yang}

\address{Mingwei Wang,  Qiuzhen College, Tsinghua University, Beijing, 100084, China}
\email{}
\address{Xiaokui Yang, Department of Mathematics and Yau Mathematical Sciences Center, Tsinghua University, Beijing, 100084, China}
\email{xkyang@mail.tsinghua.edu.cn}

\begin{document}

    \begin{abstract}  In this paper we establish  new Bochner-Kodaira  formulas with quadratic curvature terms on compact K\"ahler manifolds: for any $\eta\in \Om^{p,q}(M)$,
        $$  \la\Delta_{\bar \p} \eta,\eta\ra =\la \Delta_{\bp_F} \eta,\eta\ra +\frac{1}{4}\left\langle \left(\mathcal {R} \otimes \mathrm{Id}_{\Lambda^{p+1,q-1}T^*M}\right)(\mathbb T_\eta),\mathbb T_\eta \right\rangle.
        $$
        This linearized curvature term yields new vanishing theorems and provides estimates for Hodge numbers under exceptionally weak curvature conditions. Furthermore, we  derive Weitzenb\"ock formulas with quadratic curvature terms on both Riemannian and Kähler manifolds.
    \end{abstract}

    \maketitle

    \setcounter{tocdepth}{1}
{\small{    \begin{spacing}{1.08} \tableofcontents %
            \dottedcontents{section}[1cm]{}{2em}{5 pt} %
\end{spacing} }}

    \section{Introduction}
On a Riemannian manifold $(M,g)$, the classical Weitzenb\"ock formula states that for any differential form $\omega\in \Om^k(M)$,
    \beq  \Delta \omega=D^*D\omega+ \Theta \left(\omega\right) \label{wei}\eeq
    where $\Theta$ is a tensor induced by the curvature tensor of $(M,g)$. By combining this identity with curvature positivity conditions, numerous rigidity and vanishing theorems have been established in both Riemannian and K\"ahler geometry, e.g., \cite{Bochner46,Bochner48,Bochner49, EellsSampson1964,Meyer1971,GM75,Siu1980, SiuYau1980,Hamilton1982, Hamilton1986,MicallefMoore1988, ChenZhu2006,BrendleSchoen2008,BohmWilking2008,BrendleSchoen2009, Brendle2010,ChenTangZhu2012,Brendle2019}. For a comprehensive overview of these developments, we refer to  \cite{Wu88} and \cite{Pet16}  and the references therein. \\

    The curvature term  $\Theta$ in the Weitzenb\"ock formula is fundamentally instrumental in deriving new geometric results, particularly through its role in connecting analytic and topological properties of manifolds. In their recent work, Petersen and Wink \cite{PetersenWink2021a} established novel vanishing theorems and eigenvalue-based estimates for Betti numbers on Riemannian manifolds by exploiting positivity conditions on the curvature operator. Their framework extends to K\"ahler geometry via complexification techniques, allowing for systematic treatment of the K\"ahler curvature operator in compact settings \cite{PetersenWink2021b}. A key innovation lies in their algebraic reformulation of the curvature term into specialized representations of the Lie algebras, which provides crucial structural insights while eliminating classical singularities associated with curvature-dependent operators. For more related interesting works, we refer to \cite{Wil13}, \cite{Yang2018}, \cite{PetersenWink2022}, \cite{NienhausPetersenWink2023}, \cite{Li23}, \cite{ChoDung2023}, \cite{CGT23}, \cite{BettiolGoodman24}, \cite{DaiFu2024}, \cite{Li24}, \cite{YZ25},  \cite{Xu2025}, and their associated reference lists. \\

    The Bochner-Kodaira formula for Hermitian vector bundles $(E,h)$ over compact K\"ahler manifolds $(M,\omega)$,
    \beq \Delta_{\bp_E}=\Delta_{\p_E}+\left[\sq R^E,\Lambda_\omega\right] \eeq
     represents a fundamental refinement of the Weitzenb\"ock formula. It plays a pivotal role in establishing vanishing theorems and quantitative estimates within complex algebraic geometry.
    This relationship arises because the Bochner-Kodaira technique leverages analytic methods (particularly the  $\bp$-Laplacian machinery) to derive profound global geometric consequences for holomorphic structures on complex manifolds. Its applications extend to  metric rigidity and extension problems investigations in complex differential geometry. The
    $L^2$-estimate framework derived from the formula establishes cohomology vanishing for holomorphic sections when $E$ satisfies Nakano-positivity or dual-Nakano-positivity conditions.  In \cite{Siu1980, Siu1982}, Yum-Tung Siu combined this formula with classical Weitzenb\"ock formula to obtained new vanishing theorems and rigidity theorems by using new positivity concepts.\\

    In this paper, we establish a unified framework for Weitzenb\"ock-Bochner-Kodaira formulas with transparent curvature terms in the context of abstract Hermitian vector bundles.  These formulas integrate the analytical tools of $\bp$-Laplacians with algebraic conditions on curvature tensors, and  extend results  established  in \cite{PetersenWink2021a,PetersenWink2021b}  to  broader geometric settings.

    \subsection{Bochner-Kodaira formulas with quadratic curvature terms on K\"ahler manifolds}
 To demonstrate the geometric interpretation of  curvature terms in  Bochner-Kodaira formulas, we introduce several auxiliary operators.    Let $(E,h)$ be a Hermitian holomorphic vector bundle over a K\"ahler manifold $(M,\omega_g)$.   Let $\nabla^E$ be the Chern connection of $(E,h)$ and $R^E$ be the corresponding Chern curvature. There is an induced curvature operator $\mathfrak{R}^E: \Gamma(M,T^{1,0}M \otimes E) \> \Gamma(M,T^{1,0}M \otimes E) $ given by
    \beq \left\langle \mathfrak{R}^E\left(u^{i\alpha} \frac{\p}{\p z^i} \otimes e_\alpha \right), v^{j\beta}\frac{\p}{\p z^j} \otimes e_\beta \right\rangle = R^E_{i\bar j\alpha\bar\beta}u^{i\alpha} \bar v{}^{j\beta}. \eeq
It is easy to see that $ \mathfrak{R}^E $ is a positive operator in the sense of linear algebra if and only if $ (E,h) $ is Nakano positive.  We define the contraction operator  for $ \phi \in \Omega^{p,q}(M,E) $,
    \beq \mathbb{S}_\phi: \Gamma(M,T^{*1,0}M) \> \Omega^{p,q-1}(M,E), \ \ \ \mathbb{S}_\phi(\alpha) = I_{\alpha_\#} \phi, \eeq
where $ \alpha_\# $ is the  $(0,1)$ type dual vector of  $ \alpha \in \Gamma(M,T^{*1,0}M) $.  One can view it as
    \beq \mathbb{S}_\phi \in \Gamma(M,T^{1,0}M \otimes \Lambda^{p,q-1}T^*M\otimes E )\cong\Gamma\left(M,\left(T^{1,0}M \otimes E\right) \otimes \Lambda^{p,q-1}T^*M\right).\eeq
    As analogous to the induced curvature operator $\mathfrak R^E$,
    one can define the \emph{symmetrized curvature operator}  $ \mathcal {R}: \Gamma(M,\mathrm{Sym}^2 T^{1,0}M) \> \Gamma(M,\mathrm{Sym}^2 T^{1,0}M )$ by the relation \beq  g\left( \mathcal{R} \left(a\right),  b\right) = R_{i\bar j k \bar\ell} a^{ik}\bar b{}^{j\ell} \eeq where  $a=\sum a^{ik}\frac{\p}{\p z^i} \otimes \frac{\p}{\p z^k}$ and $b=\sum b^{j\ell}\frac{\p}{\p z^j} \otimes \frac{\p}{\p z^\ell}$ are in  $\Gamma(M,\mathrm{Sym}^2 T^{1,0}M)$ (see also \cite{CV60} and \cite{BNPSW25}) .  We say that  $(M,\omega_g)$ has \emph{positive symmetrized curvature operator} $\mathcal R$
if it is positive definite as a Hermitian bilinear form.
A straightforward computation shows that the symmetrized curvature operator of $(\C\P^n, \omega_{\mathrm{FS}})$ is $\mathcal R=2\cdot \mathrm{Id}$ which is  positive definite. On the other hand,
    if  $\mathcal {R}$ is a positive operator, then $(M,g)$ has positive holomorphic bisectional curvature. Furthermore,  when $M$ is compact,  it follows from Siu-Yau's solution to the Frankel conjecture (\cite{SiuYau1980,Mori1979}) that  $M$ is biholomorphic to $\C\P^n$.\\

For any $q\geq 1 $ and  differential form $\eta\in \Om^{p,q}(M, E)$,  we define
$$  \mathbb{T}_\eta: \Om^{1,0}(M)\times \Om^{1,0}(M)\>\Om^{p+1,q-1}(M,E) $$
by the contraction formula: \beq  \mathbb T_\eta(\alpha,\beta) = I_{\beta_\#}(\alpha \wedge \eta) + I_{\alpha_\#}(\beta \wedge \eta).\eeq
    It is obvious that  $\mathbb T_\eta(\alpha,\beta) = \mathbb T_\eta(\beta,\alpha)$ and so \beq  \mathbb T_\eta \in \Gamma\left(M,\mathrm{Sym}^2T^{1,0}M\otimes \Lambda^{p+1,q-1}T^*M\ts E\right).\eeq
    The space  $ \Gamma(M,\Sym^2 T^{1,0}M \otimes \Lambda^{p,q}T^*M\ts E) $ serves as the natural domain for expressing curvature interactions in the Bochner-Kodaira framework. We equip this bundle with a fiberwise ‌Hermitian inner product‌ by:
\beq  \left\langle u \otimes \alpha, v \otimes \beta \right\rangle_{\Sym^2T^{1,0}M\otimes \Lambda^{p,q}T^*M\ts E} := g(u,v) \langle \alpha ,\beta \rangle, \eeq
where  $ u, v \in \Gamma(M, \Sym^2 T^{1,0}M) $ and $ \alpha,\beta \in\Gamma\left(M, \Lambda^{p,q}T^*M\ts E\right)$.
 The symmetrized curvature operator $ \mathcal{R} $  admits a canonical extension to the tensor product bundle $ \Gamma(M,\Sym^2 T^{1,0}M \otimes \Lambda^{p,q}T^*M\ts E ) $  via the algebraic tensor product construction
\beq  (\mathcal {R} \otimes \text{Id}_{\Lambda^{p,q}T^*M\ts E})(u\otimes\alpha) := \mathcal{R}(u) \otimes \alpha. \eeq

    We now establish the  Bochner-Kodaira formula for the symmetrized curvature operator $\mathcal R$ and $\mathfrak{R}^E $ which reveals the intricate relationship between the Laplacians and curvature interaction terms on K\"ahler manifolds.

    \btheorem   \label{KBformulaVBcase}
    Let $(M,\omega_g)$ be a compact K\"ahler manifold and $(E,h)$ be a Hermitian holomorphic vector bundle over $M$. For any $E$-valued $(p,q)$ form $ \phi \in \Omega^{p,q}(M,E) $, one has the following Bochner-Kodaira formula:
\begin{eqnarray}  \langle\Delta_{\bar \p_E}\phi,\phi\rangle = \langle \Delta_{\bp_F}  \phi, \phi\rangle \nonumber&+& \frac{1}{4}\left\langle \big(\mathcal{R}\otimes \mathrm{Id}_{\Lambda^{p+1,q-1}T^*M\otimes E}\big)\big(\mathbb{T}_{\phi}\big),\mathbb{T}_{\phi} \right\rangle\\ &+& \left\langle \big(\mathfrak{R}^E \otimes \mathrm{Id}_{\Lambda^{p,q-1}T^*M}\big) \left(\mathbb{S}_\phi\right), \mathbb{S}_\phi \right\rangle,\label{VBKB} \end{eqnarray}
    where $ F = \Lambda^{p,q}T^*M \otimes E $  is the  \emph{complex} vector bundle  and $\bp_F$ is the $(0,1)$-part of the induced metric connection on $F$.
    \etheorem

\noindent    By convention, the curvature term $$\left\langle \big(\mathcal{R}\otimes \mathrm{Id}_{\Lambda^{p+1,q-1}T^*M\ts E}\big)\big(\mathbb{T}_{\phi}\big),\mathbb{T}_{\phi} \right\rangle$$ vanishes when $q=0$ or $p=n$. Moreover, when $q=0$, $F$ is a holomorphic vector bundle, and \beq \langle\Delta_{\bar \p_E}\phi,\phi\rangle = \langle \Delta_{\bp_F}  \phi, \phi\rangle. \eeq
\noindent When $p=n$, one obtains the formulation
    \beq \langle\Delta_{\bar \p_E}\phi,\phi\rangle = \langle \Delta_{\bp_F}  \phi, \phi\rangle  + \left\langle \big(\mathfrak{R}^E \otimes \mathrm{Id}_{\Lambda^{n,q-1}T^*M}\big) \left(\mathbb{S}_\phi\right), \mathbb{S}_\phi \right\rangle. \eeq
    \noindent When $E$ is the trivial line bundle, one obtains the following Bochner-Kodaira formula on compact K\"ahler manifolds.

    \btheorem\label{main1} Let $(M^n,\omega_g)$ be a compact K\"ahler manifold.
    For any   differential form $ \phi \in \Omega^{p,q}(M) $, the following Bochner-Kodaira formula holds
    \beq  \la\Delta_{\bar \p} \phi ,\phi \ra =\la \Delta_{\bp_F} \phi ,\phi \ra +\frac{1}{4}\left\langle \left(\mathcal {R} \otimes \mathrm{Id}_{\Lambda^{p+1,q-1}T^*M}\right)(\mathbb T_\phi ),\mathbb T_\phi  \right\rangle, \label{key0}
    \eeq
    where $F=\Lambda^{p,q}T^*M$ is the  complex vector bundle of
    $(p,q)$-forms.
    \etheorem

    \noindent   It is well-known that the curvature term \beq\left\langle \left(\mathcal {R} \otimes \mathrm{Id}_{\Lambda^{p+1,q-1}T^*M}\right)\left(\mathbb T_\phi \right),\mathbb T_\phi  \right \rangle,\label{curvatureterm}\eeq  plays a key role in Bochner-Kodaira formula applications. The principal innovation of the Bochner-Kodaira formula \eqref{key0} lies in the geometric interpretation of the curvature term, which manifests as a symmetric bilinear form in the contraction operator $\mathbb T_\phi$ derived from $\phi$.  This representation establishes the curvature term as a quadratic functional of the original $(p,q)$-form $\phi$.
Through a rigorous analysis of the operator norm relationship between $\mathbb T_\phi $ and $\phi $,  we derive estimates that bound the curvature term's magnitude by  $|\phi |$. These results are of independent interest in K\"ahler geometry.
    \btheorem
    \label{main3}
    Let  $(M,\omega_g)$ be a compact K\"ahler manifold. Suppose that $ \phi \in \Omega^{p,q}(M) $ and $ v \in \Gamma\left(M,\Sym^2 T^{*1,0}M\right) $. Then 
    \beq |\mathbb T_\phi(v)|^2 \leq \frac{4(p+1)q}{p+q} |v|^2 |\phi|^2.\eeq 
    Moreover, if  there exists some $k\geq 0$, $k\neq \frac{p+q}{2}$ and primitive $\psi\in \Omega^{p-k,q-k}(M) $ such that  $\phi=L^k\psi $,
    then we have improved estimate
    \beq |\mathbb T_\phi(v)|^2 \leq \frac{4(p-k+1)(q-k)}{p+q-2k} |v|^2 |\phi|^2.\label{key00}\eeq
    \etheorem

    \noindent The explicit tensor formula \eqref{key0} and estimate \eqref{key00} enable more precise spectral analysis of the Laplacian on K\"ahler manifolds by uncovering previously inaccessible relationships between curvature tensors and harmonic forms. Let us define a Hermitian form $ \mathbb{B}^{p,q}: \Omega^{p,q}(M) \times \Omega^{p,q}(M) \> \mathbb{C} $ by
    \beq \mathbb{B}^{p,q}(\psi,\eta) = \left\langle \big(\mathcal{R}\otimes \mathrm{Id}_{\Lambda^{p+1,q-1}T^*M}\big)\big(\mathbb{T}_{\psi}\big),\mathbb{T}_{\eta} \right\rangle. \eeq
    By using \eqref{key00} we demonstrate that the positivity of   $ \mathbb{B}^{p,q} $ follows from certain  weak postivity of $ \mathcal{R} $. Specifically, when the symmetrized curvature operator $\mathcal{R}$ is $m$-positive (i.e., the sum of its $m$ smallest eigenvalues is positive), the curvature term $ \mathbb{B}^{p,q} $ becomes positive definite for appropriate values of $p$ and $q$.

    \btheorem
    \label{positivityofBpq} Let  $(M,\omega_g)$ be a compact K\"ahler manifold with
$m$-positive symmetrized curvature operator $ \mathcal{R} $.  Then $ \mathbb{B}^{p,q} $ is positive definite in the following cases:
    \bd
    \item $ q \geq p + 2 $ and $ m \leq \frac{(n-p+1)(p+q)}{2(p+1)}$;
    \item $ q = p + 1 $,  $p \leq \frac{n}{2} $ and $ m\leq \frac{n+1}{2} $;
    \item  $ q = p + 1 $,  $ \frac{n}{2}<p<n $ and $ m\leq \frac{(n-p+1)(2p+1)}{2(p+1)} $.
    \ed
    $ \mathbb{B}^{p,q} $ is semi-positive definite in the following cases:
    \bd
    \item  $0< q \leq p  \leq \frac{n}{2} $ and $ m \leq \frac{n-p+q}{2} $;
    \item  $ 0<q \leq p $, $ \frac{n}{2}<p<n $ and $ m \leq \frac{(n-p+1)(p+q)}{2(p+1)} $.
    \ed
    Moreover, $ \mathbb{B}^{p,q}(\phi,\phi) = 0 $ if and only if $ \phi = L^q\psi $ for some $ \psi \in \Omega^{p-q,0}(M) $.
    \etheorem

    \noindent  In particular, we establish new vanishing theorems and derive refined estimates for Hodge numbers on compact K\"ahler manifolds.

\btheorem\label{main5} Let $ (M,\omega) $ be a compact K\"ahler manifold. Suppose that  the symmetrized curvature operator $ \mathcal{R} $ is $ m $-positive. Then $ H_{\bar\p}^{p,q}(M,\mathbb{C}) = 0 $ if
\bd
\item  $ q \geq p + 2 $ and $   \frac{(n-p+1)(p+q)}{2(p+1)}\geq m$; or
\item $ q = p + 1 $ and $  \frac{n+1}{2}\geq m $.
\ed
Moreover, if  $m\leq n/2$, then $ H_{\bar\p}^{p,p}(M,\mathbb{C}) = \C $ for $0\leq p\leq n$.
\etheorem
\noindent The following result is a straightforward application of Theorem \ref{main5}, which is also obtained in \cite{BNPSW25}:
    \bcorollary\label{cohomologicalspace}  Let $ (M,\omega) $ be a compact K\"ahler manifold.  If $ \mathcal{R} $ is $\left\lfloor \frac{n}{2}\right\rfloor$-positive, then $M$ has the same cohomology ring as $\C\P^n$.
    \ecorollary

\noindent  It is well-known that (e.g. \cite{CV60}) if $(M,\omega)$ is the hyperquadric in $\C\P^{n+1}$ with the induced metric, then the  symmetrized curvatrue operator $ \mathcal{R} $ has eigenvalues  $$ \lambda_1 = 2-n  \qtq{and} \lambda_2 = \cdots = \lambda_{N} = 2 $$  where $ N = \frac{n(n+1)}{2} $. In particular, $\mathcal R$ is $\left(\left\lfloor \frac{n}{2}\right\rfloor+1\right) $-positive.\\

    \noindent   By employing Theorem \ref{KBformulaVBcase} and adapting the methodology from the proof of Theorem \ref{main5}, we establish a vanishing theorem that generalizes the classical Nakano vanishing theorem for abstract vector bundles.

    \btheorem\label{main7} Let $(M,\omega_g)$ be a compact K\"ahler manifold, and $(E,h)$ be a Hermitian holomorphic vector bundle of rank $r$.  If $ (E,h)$ is Nakano positive and the symmetrized curvature operator $ \mathcal{R} $ is $ m $-positive, then the cohomology groups vanish: $$ H_{\bar\p}^{p,q}(M,E) = 0, \qtq{for} q \geq 1 $$ in the following cases:
    \bd
    \item  $ p = n $;
    \item  $ q \leq p + 1 $, $ p \leq \frac{n}{2} $ and $ m\leq \frac{n-p+q}{2} $;
    \item $(p,q)$ is not in the case of $(1)$ or $(2)$ and $ m \leq \frac{(n-p+1)(p+q)}{2(p+1)} $.
    \ed
    \etheorem

    \noindent We emphasize that Theorem \ref{main7} incorporates curvature conditions on both the base manifold $(M,\omega_g)$ and the vector bundle
    $(E,h)$. Moreover, it is clear that  in the proof of case $(1)$, the curvature condition on the base manifold is redundant.  By using Theorem \ref{KBformulaVBcase} and Theorem \ref{positivityofBpq}, one can also obtain the fact that harmonic forms remain parallel under less restrictive conditions.

\vskip 1\baselineskip
    \subsection{Weitzenb\"ock formulas with quadratic curvature terms  on  Riemannian manifolds}

    Let $(M,g)$ be a compact and oriented Riemannian manifold. The curvature operator
    $\mathfrak R: \Gamma(M,\Lambda^2TM) \to \Gamma(M,\Lambda^2TM)$
    is defined as
    \beq
    g(\mathfrak R(X\wedge Y), Z\wedge W) = R(X,Y,W,Z).\eeq
    For any   differential form $\omega\in \Om^{p}(M)$,
    $ \mathbb{T}_\omega: \Gamma(M,T^*M)\times  \Gamma(M,T^*M)\>\Om^{p}(M)$ is the operator defined
    by the contraction formula: \beq  \mathbb T_\omega(\alpha,\beta) =  \alpha \wedge I_{\beta_\#} \omega - \beta \wedge I_{\alpha_\#} \omega,\eeq
    where $ \alpha_\#, \beta_\#  $ denote the dual vector fields of $\alpha$ and $\beta$ respectively.  It is obvious that  $\mathbb T_\omega(\alpha,\beta) =- \mathbb T_\omega(\beta,\alpha)$ and so \beq  \mathbb T_\omega \in \Gamma\left(M,\Lambda^2TM\otimes \Lambda^{p}T^*M\right).\eeq
    \noindent The following formula is analogous to Theorem \ref{main1}:
    \btheorem\label{main8} Let $(M,g)$ be a compact Riemannian manifold.
    For any   differential form $ \omega \in \Omega^p(M) $, the following Weitzenb\"ock formula  holds
    \beq \left\langle\Delta_d \omega, \omega\right\rangle = \left\langle D^*D \omega, \omega\right\rangle + \left\langle \big(\mathfrak{R}\otimes \mathrm{Id}_{\Lambda^pT^*M}\big)\left(\mathbb T_\omega\right),\mathbb T_\omega \right\rangle_{\Lambda^2TM \otimes \Lambda^pT^*M },\eeq
    where $D$ is the induced connection on $\Lambda^pT^*M$.
    \etheorem
    \noindent By applying this Weitzenböck formula with explicit quadratic curvature term, one can derive estimates analogous to those in Theorem \ref{main3} and obtain applications consistent with the results of the preceding subsection. For further details,  we refer to  Section \ref{sectionWR} and \cite{PetersenWink2021a}.

\vskip 1\baselineskip
        \subsection{Weitzenb\"ock formulas with quadratic curvature terms  on  K\"ahler manifolds}
Let $(M,\omega_g) $ be a compact K\"ahler manifold. The reduced (complexified)  curvature operator $ \mathscr{R}: \Gamma(M,T^{1,0}M \otimes T^{0,1}M) \> \Gamma(M,T^{1,0}M \otimes T^{0,1}M)$ is
 defined as:
\beq \left\langle \mathscr{R}\left(\frac{\p}{\p z^i} \wedge \frac{\p}{\p \bar z{}^j}\right), \frac{\p}{\p z^\ell} \wedge \frac{\p}{\p \bar z{}^k}\right\rangle = R_{i\bar j k\bar\ell}. \eeq
For any $ \phi \in \Omega^{p,q}(M) $, $ \mathbb{Y}_\phi: \Gamma(M,T^{*1,0}M) \times \Gamma(M,T^{*0,1}M) \> \Omega^{p,q}(M) $ is the operator defined
by the contraction formula
\beq \mathbb{Y}_\phi(\alpha,\beta) = \beta \wedge I_{\alpha_\#} \phi - \alpha \wedge I_{\beta_\#} \phi, \eeq
where $ \alpha_\# \in \Gamma(M,T^{0,1}M) $ and $ \beta_\# \in \Gamma(M,T^{1,0}M) $ are  dual vectors of $ \alpha \in \Gamma(M,T^{*1,0}M) $ and $ \beta \in \Gamma(M,T^{*0,1}M) $ respectively. It is obvious that
\beq  \mathbb{Y}_\phi \in \Gamma\left(M,\left(T^{1,0}M \wedge T^{0,1}M\right) \otimes \Lambda^{p,q}T^*M\right). \eeq

\noindent The following formula is a complex analogue of Theorem \ref{main8}:
\btheorem\label{main9} Let $(M,\omega_g)$ be a compact K\"ahler manifold. For any differential form $ \phi \in \Omega^{p,q}(M) $, the following Weitzenb\"ock formula  holds
\beq \langle\Delta_d \phi, \phi\rangle = \langle D^*D\phi, \phi\rangle + \left\langle \big(\mathscr{R}\otimes \mathrm{Id}_{\Lambda^{p,q}T^*M}\big)\big(\mathbb Y_\phi\big),\mathbb Y_\phi \right\rangle_{\left(T^{1,0}M \wedge T^{0,1}M\right) \otimes \Lambda^{p,q}T^*M }, \eeq
where $ D $ is the induced connection on $ \Lambda^{p,q}T^*M $.
\etheorem

\noindent We refer to Section \ref{sectionWK} and \cite{PetersenWink2021b} for more discussions.
\vskip 1\baselineskip

\noindent\textbf{Acknowledgements}. The Weitzenb\"ock-Bochner-Kodaira formulas presented in this paper were derived in August 2021. The second author wishes to acknowledge Professor Kefeng Liu, whose encouragement was instrumental in the completion of this research.

\vskip 2\baselineskip

\section{The symmetrized curvature operator}\label{sectionP}

Let $ (M,\omega_g) $ be a compact K\"ahler manifold with $\dim M= n$. In local holomorphic coordinates $\{z^i\}$ of $M$,  for any $\phi\in \Om^{p,q}(M)$,   it can be written as
\beq    \phi=\frac{1}{p!q!}\sum_{i_1,\cdots, i_p, j_1,\cdots,   j_q}\phi_{i_1\cdots i_p \bar{j_1}\cdots \bar{j_q}}dz^{i_1} \wedge   \cdots \wedge  dz^{i_p}\wedge d\bar z{}^{j_1}\wedge \cdots\wedge d\bar z^{}{j_q},\eeq
where  $\phi_{i_1\cdots i_p \bar{j_1}\cdots \bar{j_q}}$ is skew symmetric with respect to  both $i_1,\cdots, i_p$ and $j_1,\cdots, j_q$.
\noindent The {local inner product} on $\Om^{p,q}(M)$ is defined as
\beq    \label{localinnerproduct} \la\phi,\psi\ra=\frac{1}{p!q!}g^{i_1\bar {\ell_1}}\cdots g^{i_p\bar {\ell_p}}g^{k_1\bar {j_1}}\cdots g^{k_q\bar {j_q}} \phi_{i_1\cdots i_p \bar{j_1}\cdots \bar{j_q}}\cdot \bar{\psi_{\ell_1\cdots \ell_p \bar{k_1}\cdots \bar{k_q}}} \eeq
and the {norm} on $\Om^{p,q}(M)$ is given by
\beq \|\phi\|^2=(\phi,\phi)=\int_M \la\phi,\phi\ra\frac{\omega^n}{n!}.\eeq
It is well-known that there exists a real isometry $*:\Om^{p,q}(M)\>\Om^{n-q,n-p}(M)$ such that
\beq (\phi,\psi)=\int \phi\wedge *\bar\psi.\eeq
The formal adjoint operators of $\p$ and $\bp$  are denoted by $\p^*$ and $\bp^* $ respectively. \\

\noindent
 For any  $X\in \Gamma(M,TM)$, $I_X$ is the contraction operator, i.e., \beq (I_X\eta)(\bullet)=\eta(X,\bullet), \eeq
for $\eta\in\Om^\bullet(M)$. In local  coordinates,
we also write $I_{i}$ for $I_{\frac{\p}{\p z^i}}$ and $I_{\bar j}$ for $I_{\frac{\p}{\p \bar z^j}}$. For any  $\phi \in \Omega^{p,q}(M)$, the following formulas are well-known
\beq \p \phi=dz^i \wedge \nabla_i \phi, \ \ \  \bp \phi=d\bar z{}^i \wedge \nabla_{\bar i}\phi, \label{21} \eeq
and
\beq \p^* \phi= -g^{i\bar j}I_i\nabla_{\bar j}\phi, \ \ \  \bp^* \phi = -g^{i\bar j}I_{\bar j}\nabla_i \phi, \label{22}\eeq
where $\nabla $ is the  connection on $ \Lambda^{p,q}T^*M $ induced by the Levi-Civita connection. \\

Let $(E,h)$ be a Hermitian
\emph{complex} (possibly \emph{non-holomorphic}) vector bundle over $(M,\omega_g)$. Let $\nabla^E$ be an arbitrary \emph{metric
    connection} on $(E,h)$, i.e.,   for any $s,t\in \Gamma(M,E)$, \beq dh(s,t)=h(\nabla^E s, t)+h(s,
\nabla^E t). \eeq There is a natural
decomposition \beq \nabla^E=\nabla^{'E}+\nabla{''^E} \eeq where
\beq
\nabla^{'E}:\Gamma(M,E)\>\Om^{1,0}(M,E),\ \ \ \
\nabla^{''E}:\Gamma(M,E)\>\Om^{0,1}(M,E).
\eeq
There are two induced operators
$
\p_E:\Om^{p,q}(M,E)\>\Om^{p+1,q}(M,E)$ and $\bp_E:\Om^{p,q}(M,E)\>\Om^{p,q+1}(M,E) $ given by \beq
\p_E(\phi\ts s)=\left(\p\phi\right)\ts s+(-1)^{p+q}\phi\wedge
\nabla^{'E}s,\ \ \
\bp_E(\phi\ts s)=\left(\bp\phi\right)\ts s+(-1)^{p+q}\phi\wedge
\nabla^{''E}s, \eeq where $\phi\in \Om^{p,q}(M)$ and $s\in
\Gamma(M,E)$ are local sections. The following formula is well-known \beq
\left(\p_E\bp_E+\bp_E\p_E\right)(\phi\ts s)=\phi\wedge
\left(\p_E\bp_E+\bp_E\p_E\right)s.\eeq Actually,  the operator $\p_E\bp_E+\bp_E\p_E$ is
represented by the $(1,1)$  component $R^E \in \Gamma(M,
\Lambda^{1,1}T^*M\ts E^*\ts E)$ of the curvature tensor of $(E,\nabla^E)$. The norm on
$\Om^{p,q}(M,E)$ can be defined similarly. The dual operators of $\p_E$
and $\bp_E$ are denoted by $\p^*_E$ and $\bp_E^*$
respectively. The following lemma is analogous to formulas \eqref{21} and \eqref{22}.

\blemma\label{keylemma} Let $(E,h)$ be a Hermitian complex vector bundle over a compact K\"ahler manifold $(M,\omega_g)$.
For any $\phi \in \Omega^{p,q}(M,E) $, one has
\beq \p_E \phi = dz^i \wedge \hat \nabla_i \phi, \ \ \ \bp_E \phi = d\bar z{}^j \wedge \hat \nabla_{\bar j} \phi, \eeq
and
\beq \p_E^* \phi = - g^{i\bar j} I_i \hat \nabla_{\bar j} \phi, \ \ \ \bp_E^* \phi = -g^{i\bar j} I_{\bar j} \hat \nabla_i \phi, \eeq
where $\hat \nabla = \nabla^{\Lambda^{p,q}T^*M \otimes E}$ is the induced connection on the tensor bundle $ \Lambda^{p,q}T^*M \otimes E $.
\elemma
\noindent Lemma \ref{keylemma} is essentially well-known (see, e.g. \cite[Lemma~8.9]{LY12}) and  can be verified using duality methods. It will also be frequently employed in local computations.

\blemma
\label{nablawedgecontraction} Let $(E,h)$ be a Hermitian complex vector bundle over a compact K\"ahler manifold $(M,\omega_g)$.
For any $ \phi \in \Omega^{p,q}(M,E) $ and $ \alpha \in \Omega^{1,0}(M) $, one has
\beq \label{nablawedge} \nabla^{\Lambda^{p+1,q} T^*M \otimes E}_X(\alpha \wedge \phi) =\left( \nabla_X\alpha\right) \wedge \phi + \alpha \wedge \left(\nabla^{\Lambda^{p,q} T^*M \otimes E}_X\phi\right), \eeq
and
\beq\nabla^{\Lambda^{p-1,q} T^*M \otimes E}_X \left(I_Y \phi\right) = I_{\nabla_X Y} \phi + I_Y \left(\nabla^{\Lambda^{p,q} T^*M \otimes E}_X \phi\right),  \label{nablacontraction}\eeq
where  $ X \in \Gamma(M,T_\C M) $ and $ Y \in \Gamma(M,T^{1,0}M) $.
\elemma

\bproof The formula \eqref{nablawedge} follows directly from the definition of affine connection.  For  \eqref{nablacontraction},
let $ Y_1, \cdots, Y_{p-1} \in \Gamma(M,T^{1,0}M) $ and $ X_1, \cdots, X_q \in \Gamma(M,T^{0,1}M) $.  One can see clearly that
\be && \nabla^E_X \left(\omega(Y,Y_1, \cdots, Y_{p-1},X_1, \cdots, X_q)\right) \\
& = & \left(\nabla^{\Lambda^{p,q} T^*M \otimes E}_X \omega\right)\left(Y,Y_1, \cdots, Y_{p-1},X_1, \cdots, X_q\right) \\&&+ \omega\left(\nabla_X Y,Y_1, \cdots, Y_{p-1}\right)  +  \sum_{\alpha = 1}^{p-1} \omega(Y,Y_1, \cdots, \nabla_X Y_\alpha ,\cdots Y_{p-1},X_1, \cdots, X_q)\\&& + \sum_{\beta = 1}^q \omega\left(Y,Y_1, \cdots, Y_{p-1},X_1, \cdots, \nabla_X X_\beta, \cdots, X_q\right) .\ee
On the other hand,
\be && \nabla^E_X \left(\omega\left(Y,Y_1, \cdots, Y_{p-1},X_1, \cdots, X_q\right)\right)\\
&=& \nabla^E_X \left(\left(I_Y\omega)(Y_1, \cdots, Y_{p-1},X_1, \cdots, X_q\right)\right) \\
& = &\left(\nabla^{\Lambda^{p-1,q} T^*M \otimes E}_X\left( I_Y\omega\right)\right)\left(Y_1, \cdots, Y_{p-1}\right) + \sum_{\alpha = 1}^{p-1}\left(I_Y\omega\right)\left(Y_1, \cdots, \nabla_X Y_\alpha ,\cdots Y_{p-1},X_1, \cdots, X_q\right) \\
& + &\sum_{\beta = 1}^q \left(I_Y\omega\right)\left(Y_1, \cdots, Y_{p-1},X_1, \cdots, \nabla_X X_\beta, \cdots, X_q\right) .\ee
By comparing these two expressions, we obtain \eqref{nablawedgecontraction}.
\eproof

  For any $\eta\in \Om^{p,q}(M)$, one has
  \beq \Delta_{\bp} \eta=\left(\bp\bp^*+\bp^*\bp\right)\eta,  \eeq
  and $\Delta_{\p} \eta=\left(\p\p^*+\p^*\p\right)\eta.$
On the other hand, if we set $(E,\nabla^E)=\left(\Lambda^{p,q}T^*M,\nabla^{\Lambda^{p,q}T^*M}\right)$, then $\nabla^E$ is a metric compatible connection on the  \emph{complex vector bundle} $E$ with the induced Hermitian metric. In particular, for $\eta\in \Gamma(M,E)$, one has
\beq \Delta_{\bp_E} \eta=\left(\bp_E\bp_E^*+\bp^*_E\bp_E\right)\eta= \bp^*_E\bp_E \eta.\eeq
Similarly, $\Delta_{\p_E}\eta=\left(\p_E\p_E^*+\p^*_E\p_E\right)\eta= \p^*_E\p_E \eta$.
It is well-known that $ \Delta_{\bp} \eta$ and $\Delta_{\bp_E} \eta$ are related by certain Bochner-Kodaira type formulas. \\

\noindent
For the reader's convenience, we recall the following notions.

\bd \item The symmetrized curvature operator \beq  \mathcal {R}: \Gamma(M,\mathrm{Sym}^2 T^{1,0}M) \> \Gamma(M,\mathrm{Sym}^2 T^{1,0}M ) \eeq  is defined as: for any $a=\sum a^{ik}\frac{\p}{\p z^i} \otimes \frac{\p}{\p z^k}\in \Gamma(M,\mathrm{Sym}^2 T^{1,0}M)$ with $ a^{ik} $ symmetric,
\beq \mathcal {R}\left( a\right) = a^{ik} g^{p\bar j}g^{s\bar\ell} R_{i\bar jk\bar\ell} \frac{\p}{\p z^s} \otimes \frac{\p}{\p z^p}.\eeq

\item  For any $q\geq 1 $ and $\eta\in \Om^{p,q}(M, E)$,   $\mathbb{T}_\eta: \Om^{1,0}(M)\times \Om^{1,0}(M)\>\Om^{p+1,q-1}(M,E)$ is
\beq   \mathbb T_\eta(\alpha,\beta) = I_{\beta_\#}(\alpha \wedge \eta) + I_{\alpha_\#}(\beta \wedge \eta). \eeq
It is obvious that  $\mathbb T_\eta(\alpha,\beta) = \mathbb T_\eta(\beta,\alpha)$ and so \beq  \mathbb T_\eta \in \Gamma\left(M,\mathrm{Sym}^2T^{1,0}M\otimes \Lambda^{p+1,q-1}T^*M\ts E\right).\eeq

\item The induced  inner product on the space  $ \Gamma\left(M,\Sym^2 T^{1,0}M \otimes \Lambda^{p,q}T^*M\ts E\right) $ is:\beq  \left\langle u \otimes \alpha, v \otimes \beta \right\rangle_{\Sym^2T^{1,0}M\otimes \Lambda^{p,q}T^*M\ts E} := g(u,v) \langle \alpha ,\beta \rangle, \eeq
where  $ u, v \in \Gamma(M, \Sym^2 T^{1,0}M) $ and $ \alpha,\beta \in\Gamma\left(M, \Lambda^{p,q}T^*M\ts E\right)$ are local sections.

\item  The symmetrized curvature operator $ \mathcal{R} $ is extended to $ \Gamma(M,\Sym^2 T^{1,0}M \otimes \Lambda^{p,q}T^*M\ts E ) $  as
\beq  (\mathcal {R} \otimes \text{Id}_{\Lambda^{p,q}T^*M\ts E})(u\otimes\alpha) := \mathcal{R}(u) \otimes \alpha. \eeq
where $ u \in \Gamma(M, \Sym^2 T^{1,0}M) $ and $ \alpha \in\Gamma\left(M, \Lambda^{p,q}T^*M\ts E\right)$ are local sections.

\item For any $ \phi \in \Omega^{p,q}(M,E) $,
\beq \mathbb{S}_\phi: \Gamma(M,T^{*1,0}M) \> \Omega^{p,q-1}(M,E), \ \ \ \mathbb{S}_\phi(\alpha) = I_{\alpha_\#} \phi, \eeq
where $ \alpha_\# $ is the dual vector of $ \alpha \in \Gamma(M,T^{*1,0}M) $. In particular,
\beq \mathbb{S}_\phi \in \Gamma\left(M, \left(T^{1,0}M \otimes E\right) \otimes \Lambda^{p,q-1}T^*M\right). \eeq

\item The operator $ \mathfrak{R}^E: \Gamma(M,T^{1,0}M \otimes E) \> \Gamma(M,T^{1,0}M \otimes E) $ is defined by
\beq \left\langle \mathfrak{R}^E\left(u^{i\alpha} \frac{\p}{\p z^i} \otimes e_\alpha \right), v^{j\beta}\frac{\p}{\p z^j} \otimes e_\beta \right\rangle_{T^{1,0}M \otimes E} = R^E_{i\bar j\alpha\bar\beta}u^{i\alpha} \bar v{}^{j\beta} , \eeq
for any $ u = u^{i\alpha} \frac{\p}{\p z^i} \otimes e_\alpha $ and $ v = v^{j\beta} \frac{\p}{\p z^j} \otimes e_\beta $ in $ \Gamma(M,T^{1,0}M \otimes E) $.

\ed

\bremark For a compact complex manifold $M$ of complex dimension $\geq 2$, its holomorphic tangent bundle $T^{1,0}M$ cannot be Nakano positive. This follows directly from the Nakano vanishing theorem (specifically, part $(1)$ of Theorem \ref{main7}), which implies that if $T^{1,0}M$ were Nakano positive, then $M$ is K\"ahler and we would have \begin{equation} 0 = H^{n,n-1}(M,T^{1,0}M) \cong H^{1,1}(M,\mathbb{C}), \end{equation} a manifest contradiction. This observation provides strong motivation for considering the symmetrized curvature operator $ \mathcal{R}:\Gamma(M,\mathrm{Sym}^2 T^{1,0}M) \rightarrow \Gamma(M,\mathrm{Sym}^2 T^{1,0}M). $  A natural question arising from this consideration is:
\eremark
\begin{conjecture} \label{conjecture}Let $ (M^n,\omega) $ be a compact K\"ahler manifold with $k$-positive symmetrized curvature operator $\mathcal R$. If $1\leq k\leq \left[ \frac{n}{2}\right]  $, then $M$ is biholomorphic to $\C\P^n$. 
\end{conjecture}

\vskip 1\baselineskip

\section{Bochner-Kodaira formulas with quadratic curvature terms on  K\"ahler manifolds}

\noindent In this section, we establish new Bochner-Kodaira formulas with quadratic curvature terms on compact  K\"ahler manifolds, thereby proving Theorem \ref{KBformulaVBcase} and Theorem \ref{main1}. For reader's convenience, we restate Theorem \ref{KBformulaVBcase} below:

\btheorem Let $(M,\omega_g)$ be a compact K\"ahler manifold and $(E,h)$ be a Hermitian holomorphic vector bundle over $M$. For any $ \phi \in \Omega^{p,q}(M,E) $, one has
\begin{eqnarray}  \langle\Delta_{\bar \p_E}\phi,\phi\rangle = \langle \Delta_{\bp_F}  \phi, \phi\rangle \nonumber&+& \frac{1}{4}\left\langle \big(\mathcal{R}\otimes \mathrm{Id}_{\Lambda^{p+1,q-1}T^*M\otimes E}\big)\big(\mathbb{T}_{\phi}\big),\mathbb{T}_{\phi} \right\rangle\\ &+& \left\langle \big(\mathfrak{R}^E \otimes \mathrm{Id}_{\Lambda^{p,q-1}T^*M}\big) \left(\mathbb{S}_\phi\right), \mathbb{S}_\phi \right\rangle, \label{key1}\end{eqnarray}
where $ F = \Lambda^{p,q}T^*M \otimes E $.
\etheorem

\bproof We establish curvature identities by using local representations of  $\bp_E^*$ and $\bp_E$ established in Section \ref{sectionP}.  For simplicity, we denote by $\nabla$ the induced connection on $ F = \Lambda^{p,q}T^*M \otimes E $ when no confusion arises. By Lemma \ref{keylemma}, one has
\be
\bar \p{}_E^* \bar \p_E \phi
& = & -g^{i\bar j} I_{\bar j} \nabla_i (d\bar z{}^k \wedge \nabla_{\bar k} \phi) \\
&=&   -g^{i\bar j} I_{\bar j} (d\bar z{}^k \wedge \nabla_i \nabla_{\bar k} \phi) \\
& = &   -g^{i\bar j} \nabla_i \nabla_{\bar j} \phi  + g^{i\bar j} d\bar z{}^k \wedge I_{\bar j} \nabla_i \nabla_{\bar k} \phi.
\ee
Similarly, one concludes that
\be
\bar \p_E \bar \p{}_E^* \phi
& = & -d\bar z{}^k \wedge \nabla_{\bar k} (g^{i\bar j} I_{\bar j} \nabla_i \phi) \\
& = & -d\bar z{}^k \wedge \left( \frac{\p g^{i\bar j}}{\p \bar z{}^k} I_{\bar j}\nabla_i \phi + g^{i\bar j}\nabla_{\bar k} I_{\bar j}\nabla_i \phi \right) \\
& = & -d\bar z{}^k \wedge \left(-g^{i\bar\ell}\bar\Gamma{}_{k\ell}^j I_{\bar j}\nabla_i \phi + g^{i\bar j} I_{\bar j}\nabla_{\bar k}\nabla_i \phi  + g^{i\bar j}\bar \Gamma{}_{jk}^\ell I_{\bar \ell}\nabla_i \phi \right)  \\
& = & -g^{i\bar j} d\bar z{}^k \wedge I_{\bar j} \nabla_{\bar k} \nabla_i \phi,
\ee
where the third identity follows from \eqref{nablacontraction}. On the other hand,
\be
\bar \p{}_F^* \bar \p_F \phi
& = & -g^{i\bar j} I_{\bar j} \nabla^{T^{*0,1}M \otimes F}_i \left(d\bar z{}^k \otimes \nabla^F_{\bar k} \phi\right)\\
& = & -g^{i\bar j} I_{\bar j} (d\bar z{}^k \otimes \nabla^F_i \nabla^F_{\bar k} \phi) \\
&=&   -g^{i\bar j} \nabla^F_i \nabla^F_{\bar j} \phi.
\ee
Since $ \nabla^F = \nabla^{\Lambda^{p,q}T^*M \otimes E}$, one can see clearly that
\beq  \nabla_i \nabla_{\bar j} \phi=\nabla^F_i \nabla^F_{\bar j} \phi.\eeq
Therefore, we obtain
\be \Delta_{\bar \p_E} \phi- \Delta_{\bar \p_F} \phi &=&\left(\bar \p{}_E^* \bar \p_E \phi +\bar \p{}_E \bar \p^*_E \phi\right)  -\left( \bar \p{}_F^* \bar \p_F \phi +\bar \p{}_F \bar \p_F^* \phi\right)\\&=& g^{i\bar j} d\bar z{}^k \wedge I_{\bar j}\left(\nabla_i\nabla_{\bar k} - \nabla_{\bar k}\nabla_i\right)\phi. \ee
If we write $ \phi = \phi^\alpha \otimes e_\alpha $ for a local frame $ \{e_\alpha\} $ of $ E $ and local forms $ \phi^\alpha \in \Omega^{p,q}(M) $,
\be &&
\Delta_{\bar \p_E} \phi- \Delta_{\bar \p_F} \phi \\
& = &  g^{i\bar j} d\bar z{}^k \wedge I_{\bar j}(\nabla_i\nabla_{\bar k} - \nabla_{\bar k}\nabla_i)\phi \\
& = &  g^{i\bar j} d\bar z{}^k \wedge I_{\bar j} \left((\nabla_i\nabla_{\bar k} - \nabla_{\bar k}\nabla_i)\phi^\alpha \otimes e_\alpha + \phi^\alpha \otimes (\nabla^E_i\nabla^E_{\bar k} - \nabla^E_{\bar k}\nabla^E_i)e_\alpha\right) \\
& = & \left( g^{i\bar j} d\bar z{}^k \wedge I_{\bar j} (\nabla_i\nabla_{\bar k} - \nabla_{\bar k}\nabla_i)\phi^\alpha \right) \otimes e_\alpha  + (g^{i\bar j} d\bar z{}^k \wedge I_{\bar j}\phi^\alpha) \otimes h^{\delta\bar \gamma} R^E_{i\bar k\alpha \bar\gamma} e_\delta. \ee
It is clear that
\beq \left(\nabla_i\nabla_{\bar k}-\nabla_{\bar k}\nabla_i\right) dz^p = - g^{p\bar m} R_{i\bar k \ell \bar m} dz^\ell,\eeq and \beq  \left(\nabla_i\nabla_{\bar k}-\nabla_{\bar k}\nabla_i\right) d\bar z^q = g^{m\bar q} R_{i\bar k  m\bar \ell } d\bar z{}^\ell. \eeq
Hence, for $ \phi^\alpha \in \Omega^{p,q}(M) $,
\beq \nabla_i\nabla_{\bar k}\phi^\alpha - \nabla_{\bar k}\nabla_i\phi^\alpha= - g^{\ell\bar n} R_{i\bar k m \bar n} dz^m \wedge I_\ell \phi^\alpha + g^{m \bar \ell} R_{i\bar k m \bar n} d\bar z{}^n \wedge I_{\bar \ell} \phi^\alpha. \eeq
By using K\"ahler symmetry, one has
\beq g^{i\bar j} d\bar z{}^k \wedge I_{\bar j}  \left( g^{m \bar\ell}R_{i\bar k m \bar n} d\bar z{}^n \wedge I_{\bar \ell} \phi^\alpha \right)
=  g^{i\bar j} g^{m\bar \ell} R_{i\bar k m \bar j} d\bar z{}^k \wedge I_{\bar \ell}  \phi^\alpha = g^{m \bar \ell} R_{m\bar k} d\bar z{}^k \wedge I_{\bar \ell}  \phi^\alpha.\eeq
Hence,
\be && g^{i\bar j} d\bar z{}^k \wedge I_{\bar j} (\nabla_i\nabla_{\bar k} - \nabla_{\bar k}\nabla_i)\phi^\alpha\\ & = & - g^{i\bar j} d\bar z{}^k \wedge I_{\bar j} \left( g^{\ell \bar n} R_{i\bar k m \bar n} dz^m \wedge I_\ell \phi^\alpha \right)+  g^{i \bar j} R_{i\bar k} d\bar z{}^k \wedge I_{\bar j}\phi^\alpha  \\
& = &  g^{i \bar j} g^{\ell\bar n} R_{i\bar k m\bar n} d\bar z{}^k \wedge I_{\bar j} I_\ell (dz^m \wedge \phi^\alpha) \\
& = &  g^{i \bar j} g^{\ell\bar n} R_{i\bar k m\bar n} I_\ell\left(d\bar z{}^k \wedge I_{\bar j} (dz^m \wedge \phi^\alpha)\right).   \ee
Therefore, one concludes
\be && \left \langle\Delta_{\bar \p_E} \phi- \Delta_{\bar \p_F} \phi , \phi \right\rangle \\
& = & h_{\alpha\bar \beta}\left\langle g^{i \bar j} g^{\ell\bar n} R_{i\bar k m\bar n} I_\ell (d\bar z{}^k \wedge I_{\bar j} (dz^m \wedge \phi^\alpha)), \phi^\beta \right\rangle + R^E_{i\bar k\alpha \bar\beta}\left\langle g^{i\bar j} d\bar z{}^k \wedge I_{\bar j}\phi^\alpha, \phi^\beta \right\rangle \\
& = & h_{\alpha\bar \beta} R_{i\bar k m \bar n} \left \langle g^{i \bar j} I_{\bar j}(dz^m \wedge \phi^\alpha), g^{k\bar \ell} I_{\bar \ell}(dz^n \wedge \phi^\beta) \right \rangle +R^E_{i\bar k\alpha \bar\beta}\left\langle g^{i\bar j} I_{\bar j}\phi^\alpha, g^{k\bar\ell}I_{\bar\ell} \phi^\beta \right\rangle. \ee
On the other hand,
\beq \mathbb T_\phi (dz^i, dz^j) = g^{i\bar k}I_{\bar k}(dz^j \wedge \phi) + g^{j\bar k} I_{\bar k}(dz^i \wedge \phi). \eeq
Hence, for any $ \phi \in \Omega^{p,q}(M,E) $,
\beq \mathbb T_\phi= \sum_{i,j} \left( \frac{\p}{\p z^i} \otimes \frac{\p}{\p z^j} + \frac{\p}{\p z^j} \otimes \frac{\p}{\p z^i} \right) \otimes g^{i\bar k}I_{\bar k}(dz^j \wedge \phi). \label{calculationofT} \eeq
In particular, we have
\beq \left\langle \big(\mathcal{R}\otimes \mathrm{Id}_{\Lambda^{p+1,q-1}T^*M\otimes E}\big)\big(\mathbb{T}_{\phi}\big),\mathbb{T}_{\phi} \right\rangle= 4h_{\alpha\bar\beta}R_{i\bar k j \bar \ell} \left \langle g^{i\bar m}I_{\bar m}(dz^j \wedge \phi^\alpha), g^{k\bar n}I_{\bar n}(dz^\ell \wedge \phi^\beta) \right\rangle. \eeq
On the other hand, for $ \phi \in \Omega^{p,q}(M,E)$, one has
$\mathbb{S}_\phi(dz^i) = g^{i\bar j}I_{\bar j} \phi$
and so
\beq \mathbb{S}_\phi = \frac{\p}{\p z^i} \otimes e_\alpha \otimes g^{i\bar j}I_{\bar j}\phi^\alpha. \eeq
Therefore,
\beq \left\langle \big(\mathfrak{R}^E \otimes \mathrm{Id}_{\Lambda^{p,q-1}T^*M}\big) \left(\mathbb{S}_\phi\right), \mathbb{S}_\phi \right\rangle = R^E_{i\bar k\alpha \bar \beta}\langle g^{i\bar j}I_{\bar j}\phi^\alpha, g^{k\bar\ell}I_{\bar\ell}\phi^\beta \rangle. \eeq
In conclusion, we obtain the Bochner-Kodaira formula \eqref{key1} and complete the proof of Theorem \ref{KBformulaVBcase}. \eproof

\noindent Theorem \ref{main1} represents a specific instance of Theorem \ref{KBformulaVBcase}, where the vector bundle $E$ reduces to the trivial line bundle; this particular case will be fundamental to our subsequent analysis.

    \btheorem Let $(M,\omega_g)$ be a compact K\"ahler manifold.
For any $ \eta \in \Omega^{p,q}(M) $,
\beq  \langle\Delta_{\bar \p} \eta, \eta\rangle = \langle \Delta_{\bp_F}   \eta,  \eta\rangle +\frac{1}{4}\left\langle \big(\mathcal{R}\otimes \mathrm{Id}_{\Lambda^{p+1,q-1}T^*M}\big)\big(\mathbb{T}_{ \eta}\big),\mathbb{T}_{ \eta} \right\rangle, \label{Bochner-Kodaira}
\eeq
where $F=\Lambda^{p,q}T^*M$.
\etheorem


\vskip 2\baselineskip

\section{Proof of Theorem \ref{main3}}
Let $(M,\omega)$ be a compact K\"ahler manifold. For any $\eta\in \Om^{p,q}(M)$, $L\eta=\omega\wedge \eta$ and $\Lambda$ is the dual operator of $L$ with respect to $\omega$. In the section, we prove Theorem \ref{main3}:

\btheorem
\label{normTestimate}
 Let  $(M,\omega_g)$ be a compact K\"ahler manifold. Suppose that $ \phi \in \Omega^{p,q}(M) $ and $ v \in \Gamma\left(M,\Sym^2 T^{*1,0}M\right) $. Then 
\beq |\mathbb T_\phi(v)|^2 \leq \frac{4(p+1)q}{p+q} |v|^2 |\phi|^2.\eeq 
Moreover, if  there exists some $k\geq 0$, $k\neq \frac{p+q}{2}$ and primitive $\psi\in \Omega^{p-k,q-k}(M) $ such that  $\phi=L^k\psi $,
then we have imporved estimate
\beq |\mathbb T_\phi(v)|^2 \leq \frac{4(p-k+1)(q-k)}{p+q-2k} |v|^2 |\phi|^2. \label{keyinequality2}\eeq
\etheorem

\noindent Theorem \ref{normTestimate} is purely a local statement and can be formulated at any fixed point. We need some general computations on norm of $\mathbb T_\phi$.

\blemma
For any $ \phi \in \Omega^{p,q}(M) $, one has
\beq \label{normT1}|\mathbb T_\phi|^2 = 2(q+1)(n-p)|\phi|^2 - 2|L\phi|^2. \eeq
\elemma

\bproof
By formula \eqref{calculationofT}, one has
\be  |\mathbb T_\phi|^2
& = & 2(g_{i\bar k}g_{j \bar \ell} + g_{i\bar \ell}g_{j \bar k}) \left \langle g^{i\bar m}I_{\bar m}(dz^j \wedge \phi), g^{k\bar n}I_{\bar n}(dz^\ell \wedge \phi) \right \rangle \\
& = & 2(g_{i\bar k}g_{j \bar \ell} + g_{i\bar \ell}g_{j \bar k}) \left \langle g^{n\bar \ell}I_{n}\left ( d\bar z{}^k \wedge g^{i\bar m}I_{\bar m}(dz^j \wedge \phi) \right) , \phi \right \rangle. \ee
Moreover, a straightforward calculation shows that
\be
&& (g_{i\bar k}g_{j \bar \ell} + g_{i\bar \ell}g_{j \bar k}) \cdot g^{n\bar \ell}I_{n}\left ( d\bar z{}^k \wedge g^{i\bar m}I_{\bar m}(dz^j \wedge \phi) \right) \\
& = & I_j(d\bar z{}^k \wedge I_{\bar k}(dz^j \wedge \phi)) + g_{j\bar k} g^{i\bar m} I_i(d\bar z{}^k \wedge I_{\bar m}(dz^j \wedge \phi)).
\ee
On the other hand, one can see clearly that
\beq d\bar z{}^k \wedge I_{\bar k}(dz^j \wedge \phi)=q\left(dz^j \wedge \phi\right), \eeq and \beq  I_j\left(dz^j \wedge \phi\right)=(n-p)\phi.  \eeq
Therefore,
\beq  I_j(d\bar z{}^k \wedge I_{\bar k}(dz^j \wedge \phi))=q(n-p)\phi. \label{normTa}\eeq
We also have
\be g_{j\bar k} g^{i\bar m} I_i(d\bar z{}^k \wedge I_{\bar m}(dz^j \wedge \phi))&=&- g_{j \bar k} g^{i \bar m}d\bar z{}^k \wedge dz^j\wedge  \left(I_i I_{\bar m} \phi\right) + d\bar z{}^k \wedge I_{\bar k} \phi\\
&=&- L\Lambda  \phi+q\phi.\ee
Note also that $\left(\Lambda L-L\Lambda\right) \phi=(n-p-q)\phi$ and so
\beq g_{j\bar k} g^{i\bar m} I_i(d\bar z{}^k \wedge I_{\bar m}(dz^j \wedge \phi))=(n-p)\phi-\Lambda L\phi.\label{normTb}\eeq  By \eqref{normTa} and \eqref{normTb}, we obtain
\beq|\mathbb T_\phi|^2 = 2(q+1)(n-p)|\phi|^2 - 2\la\Lambda L\phi,\phi\ra. \eeq
This is \eqref{normT1}.
\eproof
\noindent We shall deal with the term $|L\phi|^2$ in \eqref{normT1}.

\blemma\label{normT2} Suppose that at some point $ x \in M $, $ \phi \in \Omega^{p,q}(M) $ can be written as \beq  \phi = L^k \psi \qtq{and} \Lambda \psi = 0, \eeq for some  $ 0 \leq k \leq \min\{p,q\} $ and $ \psi \in \Lambda^{p-k,q-k}T^*_xM $, then
\beq |\mathbb T_\phi|^2(x) = 2(q - k)(n - p + k + 1) |\phi|^2(x).\eeq
\elemma

\bproof Since for any $\eta\in \Om^{s,t}(M)$, $ [\Lambda , L]\eta = (n - s - t) \eta$, one has
\be
\Lambda L^k\phi - L^k \Lambda\phi
& = & \sum_{i = 0}^{k-1} L^i (\Lambda L - L\Lambda) L^{k-i-1}\phi \\
& = & \sum_{i = 0}^{k-1} L^i\left( \left(n - p - q - 2(k - i - 1)\right) L^{k-i-1}\phi\right) \\
& = & k (n - p - q - k + 1) L^{k-1}\phi.
\ee
In particular, for $ \psi \in \Lambda^{p-k,q-k}T^*_xM $, one has
\beq \Lambda L^{k+1} \psi - L^{k+1} \Lambda \psi = (k+1)(n-p-q+k) L^k \psi. \eeq
Since $\phi=L^k\psi$ and $\Lambda\psi=0$, this implies
\beq \Lambda L \phi = (k+1)(n-p-q+k) \phi. \eeq
Therefore, we conclude
\be |\mathbb T_\phi|^2 &=&2\left[(q+1)(n-p)- (k+1)(n-p-q+k)\right]|\phi|^2\\&=& 2(q - k)(n - p + k + 1) |\phi|^2.\ee
This completes the proof.
\eproof

\vone

\bproof[Proof of Theorem \ref{normTestimate}] Fix $x\in M$ and we shall verify \eqref{keyinequality2} at $x$. Let $\{\tilde z^i\}$ be a local holomorphic coordinate system centered at $x\in M$ with $g(d\tilde z^i, d\tilde z^j)(x)=\delta_{ij}$. Then $v=\tilde v_{ij}d\tilde z^i\ts d\tilde z^j$ with $\tilde v_{ij}$ symmetric. By Takagi decomposition for symmetric complex matrices (e.g. \cite[Theorem~4.5.15]{HJ85}), there exists a unitary matrix $U$ and a diagonal matrix $\Lambda=\mathrm{diag}(\lambda_1,\cdots, \lambda_n)$ with $\lambda_i\geq 0$ such that
$$(\tilde v_{ij})=U\Lambda U^T.$$
If we set $z=U^T\tilde z$, then at point $x\in M$, one has
\beq v=\sum_{i=1}^n\lambda_i dz^i\ts dz^i, \eeq
and $\{dz^i\}$ is orthonormal at point $x\in M$.
In particular,
\beq |v|^2 = \sum_{i=1}^n |\lambda_{i}|^2. \eeq
By formula \eqref{calculationofT},
\beq \label{Tephi} \mathbb T_\phi(v) = 2\sum_i\lambda_{i} I_{\bar i}(dz^i \wedge \phi).\eeq
\noindent We introduce the following notation to simplify computations. Let $ \mathscr{K} $ be the collection of all ordered index pairs  $(K_1, K_2) $ satisfying
\bd \item $K_1$ and $K_2$ are ordered subsets of $\{1,\cdots, n\}$;
\item $ K_1 $ and $ K_2 $ are disjoint;

\item $ |K_1| + 2|K_2| = p + q $.\ed Moreover,  we define \beq  V_{(K_1,K_2)}^{p,q} :=
\mathrm{Span}_\mathbb{C} \{ dz^I \wedge d\bar z{}^J \wedge dz^{K_2}
\wedge d\bar z{}^{K_2} \ |\  I\cap J=\emptyset,  I\cup J = K_1 \}
\cap \Lambda^{p,q}T_x^{*}M, \eeq where $I$ and $J$ are ordered. Then
there is an orthogonal decomposition of $ \Lambda^{p,q}T_x^{*}M $:
\beq \Lambda^{p,q}T_x^{*}M= \bigoplus_{(K_1,K_2) \in \mathscr{K}}
V^{p,q}_{(K_1,K_2)}.\eeq The inequality \eqref{keyinequality2} will
be established by using this decomposition. We first show that for
any $ \psi \in V_{(K_1,\emptyset)}^{p_0,q_0} $, one has \beq
|\mathbb T_\psi(v) |^2 \leq \frac{4(p_0+1)q_0}{p_0+q_0}
|v|^2|\psi|^2, \label{reducedinequality} \eeq and
 the operator norm inequality holds:
\beq \label{estimateT} \nm{\mathbb T_\bullet(v)}_{V^{p_0,q_0}_{(K_1,\emptyset)}} \leq \frac{4(p_0+1)q_0}{p_0+q_0} |v|^2. \eeq
\noindent
Actually, in this case, we have $p_0+q_0\leq n$ and
\be \frac{1}{4} |\mathbb T_\psi(v) |^2
& = &\sum_{i,j} \lambda_{i} \lambda_{j} \langle I_{\bar i}(dz^i \wedge \psi), I_{\bar j}(dz^j \wedge \psi) \rangle \\
& = &\sum_{i,j} \lambda_{i} \lambda_{j} \langle d\bar z{}^j \wedge I_{\bar i} \psi, I_i(dz^j \wedge \psi) \rangle \\
& = &\sum_{i,j} \lambda_{i}^2 \langle d\bar z{}^i \wedge I_{\bar i} \phi, \phi \rangle - \lambda_{i} \lambda_{j} \langle d\bar z{}^j \wedge I_{\bar i} \psi, dz^j \wedge I_i \psi \rangle \\
& = &\sum_{i,j} \lambda_{i}^2 \langle I_{\bar i} \psi, I_{\bar i} \psi \rangle + \lambda_{i} \lambda_{j} \langle I_j I_{\bar i} \psi, I_{\bar j}I_i \psi \rangle.
\ee
For fixed $i$ and $j$, there are  elementary inequalities
\beq 2| \lambda_{i} \lambda_{j} \langle I_i I_{\bar j} \psi, I_{\bar i} I_j \psi \rangle | \leq  \lambda_{i}^2 \langle I_j I_{\bar i} \psi, I_j I_{\bar i} \psi \rangle + \lambda_{j}^2 \langle I_i I_{\bar j} \psi, I_i I_{\bar j} \psi \rangle, \eeq
and
\beq 2|\lambda_{i} \lambda_{j} \langle I_i I_{\bar j} \psi, I_{\bar i} I_j \psi \rangle | \leq  \lambda_{j}^2 \langle I_j I_{\bar i} \psi, I_j I_{\bar i} \psi \rangle + \lambda_{i}^2 \langle I_i I_{\bar j} \psi, I_i I_{\bar j} \psi \rangle. \eeq
Note also that for any $\eta\in \Lambda^{s,t}T^*_xM$, $\sum_j\la I_j \eta, I_j \eta\ra=s|\eta|^2$ and $\sum_j\la I_{\bar j} \eta, I_{\bar j} \eta\ra=t|\eta|^2$. Hence,
\be \sum_{i,j} 2|\lambda_{i} \lambda_{j} \langle I_i I_{\bar j} \psi, I_{\bar i} I_j \psi \rangle | &\leq&  p_0\sum_i\lambda_{i}^2 \langle I_{\bar i} \psi, I_{\bar i} \psi \rangle + p_0\sum_j\lambda_{j}^2 \langle I_{\bar j} \psi, I_{\bar j} \psi \rangle\\&=&2p_0\sum_i\lambda_{i}^2 \langle I_{\bar i} \psi, I_{\bar i} \psi \rangle, \ee
and
\beq \sum_{ij} 2|\lambda_{i} \lambda_{j} \langle I_i I_{\bar j} \psi, I_{\bar i} I_j \psi \rangle |  \leq 2\sum_j q_0\lambda_{j}^2 \langle I_j \psi, I_j \psi \rangle.\eeq
An average argument shows
\beq\sum_{ij} |\lambda_{i} \lambda_{j} \langle I_i I_{\bar j} \psi, I_{\bar i} I_j \psi \rangle | \leq p_0\cdot \frac{q_0-1}{p_0+q_0}\sum_i \lambda_{i}^2 \langle I_{\bar i} \psi, I_{\bar i} \psi \rangle +q_0\cdot  \frac{p_0+1}{p_0+q_0}\sum_i \lambda_{i}^2 \langle I_i \psi, I_i \psi \rangle. \eeq
Therefore,
\be  \frac{1}{4} |\mathbb T_\psi(v) |^2
& \leq & \frac{q_0(p_0+1)}{p_0+q_0} \sum_i\left( \lambda_{i}^2 \langle I_{\bar i} \psi, I_{\bar i} \psi \rangle + \lambda_{i}^2 \langle I_i \psi, I_i \psi \rangle \right) \\
& = & \frac{q_0(p_0+1)}{p_0+q_0} \sum_i \langle \lambda_{i}^2 (d\bar z{}^i \wedge I_{\bar i} \psi + dz^i \wedge I_i \psi), \psi \rangle. \ee
Since  $ \psi \in V^{p_0,q_0}_{(K_1,\emptyset)} $ and $|K_1|\leq n$,  a straightforward computation shows that
\beq \sum_i \langle \lambda_{i}^2 (d\bar z{}^i \wedge I_{\bar i} \psi + dz^i \wedge I_i \psi), \psi \rangle= \sum_{i \in K_1} \lambda_{i}^2 |\psi|^2\leq |v|^2|\phi|^2. \eeq
Therefore, we obtain \eqref{reducedinequality}. \\

\noindent For a general pair $  (K_1,K_2) $,  we set $ (p_0,q_0) = (p,q) - (t,t) $ for $ t = |K_2| \geq 0 $. We need the following fact.\\

\noindent \emph{Claim.} There are  identities on operator norms:
\beq \nm{\mathbb T_\bullet(v)}_{V^{p_0,q_0}_{(K_1,\emptyset)}}=\nm{\mathbb T_\bullet(v)}_{V^{p,q}_{(K_1,K_2)}},\ \  \nm{\mathbb T_\bullet(v)}_{\Lambda^{p,q}T_x^{*}M}= \sup_{(K_1,K_2) \in \mathscr{K}} \nm{\mathbb T_\bullet(v)}_{V^{p,q}_{(K_1,K_2)}}. \label{normidentity}\eeq

\noindent \emph{Proof of Claim.} The linear map $ f : V^{p_0,q_0}_{(K_1,\emptyset)} \> V^{p,q}_{(K_1,K_2)}$ defined by
$ f(\alpha) = \alpha \wedge dz^{K_2} \wedge d\bar z{}^{K_2}$
is an isomorphism and
$$ f\left(\mathbb T_{\alpha}(v)\right)=\mathbb T_{f(\alpha)}(v).$$

\noindent   Indeed, it is easy to see that $f$ is surjective and $ \dim V_{(K_1,\emptyset)}^{p_0,q_0} = \dim V_{(K_1,K_2)}^{p,q} = {p_0+q_0 \choose p_0} $.  Hence, $f$ is an isomorphism. Moreover, suppose that $ K_2 = \{ k_1, \cdots, k_t \} $.
Since $ \alpha \in V_{(K_1,\emptyset)}^{p_0,q_0} $ and $ K_1 \cap K_2 = \emptyset $, for any $ i \in K_2 $, $ I_{\bar i}\alpha = 0 $ and $ I_i \alpha = 0 $. Therefore,
\be |f(\alpha)|^2
& = & \langle I_{\bar{k_t}} \cdots I_{\bar{k_1}} I_{k_t} \cdots I_{k_1} \left(dz^{k_1} \wedge \cdots \wedge dz^{k_t} \wedge d\bar z{}^{k_1} \wedge \cdots \wedge d\bar z{}^{k_t} \wedge \alpha \right), \alpha \rangle \\
& = & \langle \alpha, \alpha \rangle = |\alpha|^2. \ee

\noindent For any $ \alpha \in V_{(K_1,\emptyset)}^{p_0,q_0} $,
\be \mathbb{T}_{f(\alpha)}(v) & = & 2\sum_i \lambda_iI_{\bar i}(dz^i \wedge \alpha \wedge dz^{K_2} \wedge d\bar z{}^{K_2}) \\
& = & 2\sum_{i \notin K_2} \lambda_i I_{\bar i}(dz^i \wedge \alpha \wedge dz^{K_2} \wedge d\bar z{}^{K_2}) \\
& = & -2 \sum_{i \notin K_2} \lambda_i dz^i \wedge I_{\bar i} \alpha \wedge dz^{K_2} \wedge d\bar z{}^{K_2}. \ee
Since $ \alpha \in V_{(K_1,\emptyset)}^{p_0,q_0} $ and $ K_1 \cap K_2 = \emptyset $, for any $ i \in K_2 $,  $ I_{\bar i}\alpha = 0 $, one has
\beq -2 \sum_{i \notin K_2} \lambda_i dz^i \wedge I_{\bar i} \alpha  = -2 \sum_i \lambda_i dz^i \wedge I_{\bar i} \alpha = \mathbb{T}_\alpha(v). \eeq
Hence,  $ f(\mathbb{T}_\alpha(v)) = \mathbb{T}_{f(\alpha)}(v) $. \\

\noindent For any $ \alpha \in  V_{(K_1,\emptyset)}^{p_0,q_0} $,  $|f(\alpha)| = |\alpha|$ and
$|\mathbb{T}_\alpha(v)| = |f(\mathbb{T}_\alpha(v))| = |\mathbb{T}_{f(\alpha)}(v)|$.
Since $ f $ is an isomorphism, one has
\beq \nm{\mathbb{T}_\bullet(v)}_{V^{p_0,q_0}_{(K_1,\emptyset)}} = \sup_{0 \neq \alpha \in V^{p_0,q_0}_{(K_1,\emptyset)}} \frac{|\mathbb{T}_\alpha(v)|}{|\alpha|} = \sup_{0 \neq f(\alpha) \in V^{p,q}_{(K_1,K_2)}} \frac{|\mathbb{T}_{f(\alpha)}(v)|}{|f(\alpha)|} = \nm{\mathbb{T}_\bullet(v)}_{V_{(K_1,K_2)}^{p,q}}. \eeq
Moreover,  it is clear that
\beq \sup_{(K_1,K_2) \in \mathscr{K}} \nm{\mathbb{T}_\bullet(v)}_{V_{(K_1,K_2)}^{p,q}} \leq \nm{\mathbb{T}_\bullet(v)}_{\Lambda^{p,q}T_x^*M}. \label{norma}\eeq
We shall prove the converse. For any $ \alpha \in V_{(K_1,K_2)}^{p,q} $, then $ \mathbb{T}_\alpha(v) \in V_{(K_1,K_2)}^{p+1,q-1} $. Indeed,  without loss of generality,  we assume that $ \alpha =  dz^I \wedge d\bar z{}^J \wedge dz^{K_2} \wedge d\bar z{}^{K_2} $, where $ I \cap J = \emptyset, I \cup J = K_1, |I| = p_0, |J| = q_0 $. Since $ \mathbb{T}_\alpha(v) = 2\sum_i \lambda_i I_{\bar i}(dz^i \wedge \alpha) $, one has
\beq \mathbb{T}_\alpha(v) = 2(-1)^{p_0+1} \sum_{j \in J} \lambda_j dz^j \wedge dz^I \wedge I_{\bar j}(d\bar z{}^J) \wedge dz^{K_2} \wedge d\bar z{}^{K_2} \in V_{(K_1,K_2)}^{p+1,q-1}. \eeq
Since there are orthonormal decompositions:
\beq \Lambda^{p,q}T^*_xM = \bigoplus_{(K_1,K_2) \in \mathscr{K}} V_{(K_1,K_2)}^{p,q}, \ \ \  \Lambda^{p+1,q-1}T^*_xM = \bigoplus_{(K_1,K_2) \in \mathscr{K}} V_{(K_1,K_2)}^{p+1,q-1}, \eeq
for any $ \alpha \in \Lambda^{p,q}T^*_xM $, it can be written as
\beq \alpha = \sum_{(K_1,K_2)\in \mathscr{K}} \alpha_{(K_1,K_2)}^{p,q}, \eeq
where $ \alpha_{(K_1,K_2)}^{p,q} \in V_{(K_1,K_2)}^{p,q} $, and
\beq \mathbb{T}_\alpha(v) = \sum_{(K_1,K_2)\in \mathscr{K}} \mathbb{T}_{\alpha_{(K_1,K_2)}^{p,q}}(v), \eeq
where $ \mathbb{T}_{\alpha_{(K_1,K_2)}^{p,q}}(v) \in V_{(K_1,K_2)}^{p+1,q-1} $. Since
\beq |\alpha|^2 = \sum_{(K_1,K_2)\in \mathscr{K}} |\alpha_{(K_1,K_2)}^{p,q}|^2, \ \ \  |\mathbb{T}_\alpha(v)|^2 = \sum_{(K_1,K_2)\in \mathscr{K}} |\mathbb{T}_{\alpha_{(K_1,K_2)}^{p,q}}(v)|^2, \eeq
and
\beq |\mathbb{T}_{\alpha_{(K_1,K_2)}^{p,q}}(v) |^2 \leq \nm{\mathbb{T}_\bullet(v)}_{V_{(K_1,K_2)}^{p,q}}^2 |\alpha_{(K_1,K_2)}^{p,q}|^2, \eeq
one concludes that
\beq |\mathbb{T}_\alpha(v)|^2 \leq \sup_{(K_1,K_2) \in \mathscr{K}} \nm{\mathbb{T}_\bullet(v)}_{V_{(K_1,K_2)}^{p,q}}|\alpha|^2.\label{normb} \eeq
The equality \eqref{normidentity} follows form \eqref{norma} and \eqref{normb}. This completes the proof of Claim.\\

\noindent
By  \eqref{normidentity}, one has
\be \nm{\mathbb T_\bullet(v)}_{\Lambda^{p,q}T_x^{*}M}=\sup_{(K_1,K_2) \in \mathscr{K}} \nm{\mathbb T_\bullet(v)}_{V^{p,q}_{(K_1,K_2)}}
&=& \sup_{(K_1,K_2) \in \mathscr{K}} \nm{\mathbb T_\bullet(v)}_{V^{p_0,q_0}_{(K_1,\emptyset)}}\\
& \leq & \sup_t \left( \frac{4(p+1-t)(q-t)}{p+q-2t} \right)|v| \\
& = & \frac{4(p+1)q}{p+q}|v| . \ee
That is, for any  $ \phi \in \Omega^{p,q}(M) $ and $ v \in \Gamma(M,\Sym^2 T^{*1,0}M) $, one has
\beq |\mathbb T_\phi(v)|^2 \leq \frac{4(p+1)q}{p+q} |v|^2 |\phi|^2. \label{keyinequality1}\eeq

For any $\eta\in\Om^{a, b}(M) $, it is straightforward to verify that
\beq L\left(\mathbb T_{\eta}(v)\right)=\mathbb T_{L(\eta)}(v), \qtq{and}  \Lambda\left(\mathbb T_{\eta}(v)\right)=\mathbb T_{\Lambda(\eta)}(v). \eeq
Moreover, if  $ \Lambda \eta = 0 $, then by using a similar argument as in the proof of Lemma \ref{normT2} (e.g. Lemma \ref{LLambdaLemma}), one has
\beq \Lambda^kL^k \eta= \prod_{i=1}^k i(n+1-a-b-i) \eta \qtq{and} |L^k\eta|^2 = \prod_{i=1}^k i(n+1-a-b-i) |\eta|^2.\label{key3} \eeq

\noindent
If $\phi=L^k\psi$ and $\Lambda\psi=0$ for some $k\neq \frac{p+q}{2}$, then $\Lambda\left(\mathbb T_{\psi}(v)\right)=\mathbb T_{\Lambda(\psi)}(v)=0$ and
\beq |\mathbb T_{L^k(\psi)}(v)|^2=|L^k\left(\mathbb T_{\psi}(v)\right)|^2= \prod_{i=1}^k i(n+1-(p-k+1)-(q-k-1)-i)|\left(\mathbb T_{\psi}(v)\right)|^2.\eeq
By inequality \eqref{keyinequality1},
\beq|\left(\mathbb T_{\psi}(v)\right)|^2 \leq \frac{4(p-k+1)(q-k)}{p+q-2k} |v|^2 |\psi|^2,\eeq
and so\be |\mathbb T_{L^k(\psi)}(v)|^2
&\leq &  \prod_{i=1}^k i\left(n+1-(p-k)-(q-k)-i\right)\frac{4(p-k+1)(q-k)}{p+q-2k} |v|^2 |\psi|^2\\
&=& \frac{4(p-k+1)(q-k)}{p+q-2k} |v|^2|L^k\psi|^2,\ee
where the last inequality follows from \eqref{key3}. Hence, we establish \eqref{keyinequality2}.
\eproof


\vskip 2\baselineskip

\section{Proofs of Theorem \ref{positivityofBpq}, Theorem \ref{main5} and Theorem \ref{main7}}
\noindent In this section, we prove Theorem \ref{positivityofBpq}, Theorem \ref{main5}, Corollary \ref{cohomologicalspace} and Theorem \ref{main7}. Let us recall the $m$-positivity.
\bdefinition Let $A$ be a Hermitian $ n \times n $ matrix and $ \lambda_1 \leq \cdots \leq \lambda_n $ be eigenvalues of $ A $.   It  is said to be $ m $-positive if
\beq \lambda_1 + \cdots + \lambda_m > 0. \eeq
The symmetrized curvature operator $ \mathcal{R}: \Gamma(M,\Sym^2 T^{1,0}M) \> \Gamma(M,\Sym^2 T^{1,0}M) $ is called $ m $-positive if $ \mathcal{R} $ is $m $-positive at every point  of $ M $. One can define $m$-semi-positivity, $m$-negativity and $m$-semi-negativity in similar ways.
\edefinition
\noindent Let $A$ be an $m$-positive Hermitian $ n \times n $ matrix.  Suppose that
 $ \{e_i\}_{i=1}^n $ is an orthonormal frame of $ \C^{n} $, then
\beq  \sum_{s=1}^k \langle Ae_{i_s}, e_{i_s} \rangle \geq \lambda_1 + \cdots + \lambda_k, \label{partialtrace} \eeq
for any $1\leq i_1<\cdots<i_k\leq n$.
\vone

\noindent We need the following technical result.

\blemma\label{combinatorial2} Assume that $ 0 \leq p < n $ and $ 0 < q \leq n $. Let $ s = \min\{ p,q-1 \} $. For $ 0 \leq k \leq s $, we define
\beq C_{p,q}^k = \frac{(n-p+k+1)(p+q-2k)}{2(p+1-k)}.\eeq
\bd
\item If $ q \geq p + 2 $ or $ p > n/2 $, one has
\beq \min_{0 \leq k \leq s} C_{p,q}^k =C_{p,q}^0= \frac{(n-p+1)(p+q)}{2(p+1)}. \eeq
\item If $ q \leq p + 1 $ and $ p \leq n/2 $, one has
\beq \min_{0 \leq k \leq s} C_{p,q}^k = C_{p,q}^{q-1}=\frac{n-p+q}{2}. \eeq
\ed

\elemma

\bproof
$(1)$. Let $ m = p+1-k > 0 $.  One has
\be C_{p,q}^k & = & \frac{1}{2m}(n+2-m)(q-p-2+2m) \\
& = & -m + \frac{(n+2)(q-p-2)}{2m} + n+2 + \frac{p+2-q}{2}. \ee
If $ q \geq p+2 $, then $ C_{p,q}^k $ is decreasing for $ m $ and increasing for $ k $, and so
\beq \min_{0 \leq k \leq s} C_{p,q}^k = C_{p,q}^0 = \frac{(n-p+1)(p+q)}{2(p+1)}. \eeq
Let us consider the case $p>n/2$. In this case, if $q\geq p+2$, we are done. If   $ q \leq p + 1 $, one has $s=q-1$. Moreover,
\beq C_{p,q}^k = -m + \frac{(n+2)(q-p-2)}{2m} + n+2 + \frac{p+2-q}{2}. \eeq
It is obvious that $ C_{p,q}^k $ is convex in $ m > 0 $, and it attains its minimum when $ m $ attains the boundary. In particular, one has
\beq \min_{0\leq k\leq s}  C_{p,q}^k = \min\left\{C_{p,q}^0, C_{p,q}^{q-1}\right\} = \min\left\{ \frac{(n-p+1)(p+q)}{2(p+1)}, \frac{(n-p+q)}{2} \right\}.\eeq
Since $ q \geq 1 $ and $ p > n/2 $, one can see that
\beq \frac{(n-p+1)(p+q)}{2(p+1)} \leq \frac{(n-p+q)}{2}, \eeq
and so
\beq \min_{0 \leq k \leq s} C_{p,q}^k = C_{p,q}^{0}= \frac{(n-p+1)(p+q)}{2(p+1)}. \eeq

$(2)$. If $ q \leq p + 1 $ and $ p \leq n/2 $, one has
\beq \min_{0\leq k\leq s}  C_{p,q}^k = \min\left\{C_{p,q}^0, C_{p,q}^{q-1}\right\} = \min\left\{ \frac{(n-p+1)(p+q)}{2(p+1)}, \frac{(n-p+q)}{2} \right\}.\eeq
Since $ q \geq 1 $ and $ p \leq n/2 $, one has
\beq \frac{(n-p+1)(p+q)}{2(p+1)} \geq \frac{(n-p+q)}{2}, \eeq
and therefore
\beq \min_{0 \leq k \leq s} C_{p,q}^k = C_{p,q}^{q-1} = \frac{n-p+q}{2}. \eeq
This completes the proof.
\eproof

\blemma
\label{LLambdaLemma}
Suppose that $ \psi \in \Lambda^{p,q} T_x^*M $ and $ \Lambda \psi = 0 $. Then
\beq \Lambda^k L^k \psi = c_k \psi,  \eeq
where   $c_k= c(p,q,n,k) = \prod_{i=1}^k i(n-p-q-i+1)$ is a constant. 
Moreover, \bd \item  if $ p + q > n $, then $ \psi = 0 $.
\item  if $ p + q \leq n $, then $ c_k \geq 0 $, and  $ c_k = 0 $ if and only if $ k \geq n - p - q +1 $.\ed 
\elemma

\bproof We prove it by induction.  For  $ k = 1 $, since $ [\Lambda, L]\psi = (n-p-q)\psi $, one has 
\beq \label{identityLambdaL} \Lambda L\psi = (n-p-q)\psi.  \eeq 
Suppose that the result holds for all $t$ satisfying $ 1 \leq t < k $.  Since we established in the proof of Lemma \ref{normT2} that 
\beq
\Lambda L^k\psi - L^k \Lambda\psi
= k (n - p - q - k + 1) L^{k-1}\psi, 
\eeq
and $ \psi $ is primative,  we conclude that 
\beq \Lambda L^k \psi = k(n-p-q-k+1)L^{k-1}\psi. \eeq
By induction, one has 
\beq \Lambda^{k-1}L^{k-1} \psi = \prod_{i=1}^{k-1} i(n-p-q-i+1) \cdot \psi,  \eeq
and so
\be \Lambda^kL^k \psi &=&\Lambda^{k-1}\left(\Lambda L^k \psi\right)\\&=& k(n-p-q-k+1)\Lambda^{k-1}L^{k-1} \psi\\ &=& \prod_{i=1}^k i(n-p-q-i+1) \cdot \psi. \ee
It is easy to see that  if $ p + q \leq n $, then $ c_k \geq 0 $, and  $ c_k = 0 $ if and only if $ k \geq n - p - q +1 $. Moreover, if $ p + q  > n $, by \eqref{identityLambdaL}, one has 
\beq |L\psi|^2 + (p+q-n)|\psi|^2 = 0, \eeq
and it implies $ \psi = 0 $. 
\eproof

\noindent Theorem \ref{positivityofBpq} states that:

    \btheorem Let  $(M,\omega_g)$ be a compact K\"ahler manifold with
$m$-positive symmetrized curvature operator $ \mathcal{R} $.  Then $ \mathbb{B}^{p,q} $ is positive definite in the following cases:
\bd
\item $ q \geq p + 2 $ and $ m \leq \frac{(n-p+1)(p+q)}{2(p+1)}$;
\item $ q = p + 1 $,  $p \leq \frac{n}{2} $ and $ m\leq \frac{n+1}{2} $;
\item  $ q = p + 1 $,  $ \frac{n}{2}<p<n $ and $ m\leq \frac{(n-p+1)(2p+1)}{2(p+1)} $.
\ed
$ \mathbb{B}^{p,q} $ is semi-positive definite in the following cases:
\bd
\item  $0< q \leq p  \leq \frac{n}{2} $ and $ m \leq \frac{n-p+q}{2} $;
\item  $ 0<q \leq p $, $ \frac{n}{2}<p<n $ and $ m \leq \frac{(n-p+1)(p+q)}{2(p+1)} $.
\ed
Moreover, $ \mathbb{B}^{p,q}(\Phi,\Phi) = 0 $ if and only if $ \Phi = L^q\psi $ for some $ \psi \in \Omega^{p-q,0}(M) $.
\etheorem
\bproof  For any $ x \in M $, since $ \mathcal{R} $ is Hermitian, there exists an orthonormal frame $ \{ e_A \}_{A=1}^N  $ of $ (\Sym^2 T^{1,0}_xM, h_x ) $ such that
\beq \mathcal {R}(e_A) = \lambda_A e_A, \eeq
where $ N = n(n+1)/2 $ and $ \lambda_1 \leq \cdots \leq \lambda_N $. Let $ \{ e^A \} $ be the dual frame of $ \{ e_A \} $, then for any  $ \phi \in \Omega^{p.q}(M) $
\beq \label{localT} \mathbb T_\phi = e_A \otimes \mathbb T_\phi(e^A). \eeq
In particular, for any $ \psi,\eta \in \Omega^{p,q}(M) $, one has
\beq \mathbb{B}^{p,q}(\psi,\eta) = \left\langle \big(\mathcal{R}\otimes \mathrm{Id}_{\Lambda^{p+1,q-1}T^*M}\big)\big(\mathbb{T}_{\psi}\big),\mathbb{T}_{\eta} \right\rangle= \sum_A \lambda_A \left\langle \mathbb{T}_\psi(e^A), \mathbb{T}_\eta(e^A) \right\rangle. \eeq
For any $\eta\in\Om^{a, b}(M) $, it is straightforward to verify that
\beq \label{communicateLandT} L\left(\mathbb T_{\eta}(v)\right)=\mathbb T_{L(\eta)}(v), \qtq{and}  \Lambda\left(\mathbb T_{\eta}(v)\right)=\mathbb T_{\Lambda(\eta)}(v). \eeq
In particular, for  $ \psi \in \Omega^{p-1,q-1}(M) $ and $ \phi \in \Omega^{p,q}(M) $, one has
\be \mathbb{B}^{p,q}(L(\psi),\phi) & = & \sum_A \lambda_A \left\langle \mathbb{T}_{L(\psi)}(e^A), \mathbb{T}_\phi(e^A) \right\rangle = \sum_A \lambda_A \left\langle L \left(\mathbb{T}_\psi(e^A)\right), \mathbb{T}_\phi(e^A) \right \rangle \\
& = & \sum_A \lambda_A \left\langle \mathbb{T}_\psi(e^A), \Lambda\left(\mathbb{T}_\phi(e^A)\right) \right\rangle = \sum_A \lambda_A \left\langle \mathbb{T}_\psi(e^A), \mathbb{T}_{\Lambda(\phi)}(e^A) \right \rangle \\
& = & \mathbb{B}^{p-1,q-1}(\psi,\Lambda(\phi)). \ee \noindent For
any $ \Phi \in \Lambda^{p,q} T^*_xM $, by Lefschetz decomposition for inner product vector spaces (e.g. \cite[Proposition~1.2.30]{Huy05}),
one has \beq \Phi=\sum_{k=0}^t L^k\psi_k, \label{LD}\eeq where $ t = \min\{p,q\}
$ and $\psi_k\in  \Lambda^{p-k,q-k} T^*_xM  $ are primitive. Since $
\Lambda\psi_k = 0 $, one has \beq \Lambda^kL^k \psi_k = c_k \psi_k,
\label{LAk}\eeq where $ c_k =c(p-k,q-k,n,k)$ are constants given in Lemma \ref{LLambdaLemma}.   Moreover, if $ p+q-2k > n $, then $$ \psi_k = 0. $$ If $ p+q -2k \leq n $, then $ c_k \geq 0 $. Moreover,  if $c_k=0$, we have $\Lambda^kL^k\psi_k=0$ and so $L^k\psi_k=0$. In this case, we can choose $\psi_k=0$ in the decomposition \eqref{LD}. In the following computations, we only consider those $\psi_k$  that satisfy $c_k>0$.\\

Since $\psi_k$ is primitive, for any $
\ell > k $, by \eqref{LAk},  one has $ \Lambda^\ell L^k \psi_k = 0 $ and \beq
\mathbb{B}^{p,q}(L^k \psi_k, L^\ell \psi_\ell) =
\mathbb{B}^{p-\ell,q-\ell}(\Lambda^\ell L^k \psi_k, \psi_\ell) = 0.
\eeq Therefore one can conclude that \beq \label{computationB} \mathbb{B}^{p,q}(\Phi,\Phi) = \sum_{k} c_k \mathbb{B}^{p-k,q-k}(\psi_k,\psi_k) = \sum_{c_k>0} c_k \mathbb{B}^{p-k,q-k}(\psi_k,\psi_k). \eeq By
Lemma \ref{normT2}, one has \beq \label{sumofphi} \sum_A |\mathbb
T_{\psi_k}(e^A)|^2 = |\mathbb T_{\psi_k}|^2 = 2(q-k)(n-p+k+1)
|\psi_k|^2. \eeq  Moreover, by \eqref{keyinequality2}, one obtains \beq
\label{norminequation} |\mathbb T_{\psi_k}(e^A)|^2 \leq
\frac{4(q-k)(p+1-k)}{p+q-2k}|e^A|^2 |\psi_k|^2. \eeq  On the other hand, since
$ \lambda_N\geq \cdots \geq \lambda_{m+1}>0$  at $ x \in M $, one
has \beq \mathbb{B}^{p-k,q-k}(\psi_k,\psi_k) =\sum_{A} \lambda_A |\mathbb T_{\psi_k}(e^A)|^2 \geq
\sum_{i=1}^m \lambda_i|\mathbb T_{\psi_k}(e^i)|^2 +
\lambda_{m+1}\sum_{j=m+1}^N |\mathbb T_{\psi_k}(e^j)|^2. \eeq By
using \eqref{sumofphi}, \beq \sum_{A} \lambda_A |\mathbb
T_{\psi_k}(e^A)|^2 \geq  \sum_{i=1}^m
(\lambda_i-\lambda_{m+1})|\mathbb T_{\psi_k}(e^i)|^2+ 2\lambda_{m+1}
(q-k)(n-p+k+1) |\psi_k|^2.\eeq Moreover, the inequality
\eqref{norminequation} gives \be \sum_{A} \lambda_A |\mathbb
T_{\psi_k}(e^A)|^2 &\geq&  4\sum_{i=1}^m
(\lambda_i-\lambda_{m+1})\frac{(q-k)(p+1-k)}{p+q-2k}|\psi_k|^2 \\&&+
2\lambda_{m+1} (q-k)(n-p+k+1)|\psi_k|^2.\ee Therefore, \beq \sum_A
\lambda_A|\mathbb T_{\psi_k}(e^A)|^2 \geq
\frac{4(q-k)(p+1-k)}{p+q-2k}\left(\left(C_{p,q}^k-m\right)\lambda_{m+1} +
\sum_{i=1}^m \lambda_i \right)|\psi_k|^2, \eeq
where $C_{p,q}^k$ is the number defined in Lemma \ref{combinatorial2}:
\beq C_{p,q}^k = \frac{(n-p+k+1)(p+q-2k)}{2(p+1-k)}.\eeq

\noindent 
 Therefore, we
conclude that if  $ \mathcal{R} $ is $ m
$-positive and $ m \leq C_{p,q}^k $, then \beq \label{PostivityB}
\mathbb{B}^{p-k,q-k}(\psi_k,\psi_k) \geq
\frac{4(q-k)(p+1-k)}{p+q-2k}\left(
\sum_{i=1}^m \lambda_i \right)|\psi_k|^2\geq 0. \eeq Moreover, if  $
\psi_k \neq 0 $ and $ k < q $, then $ \mathbb{B}^{p-k,q-k}(\psi_k,\psi_k) > 0 $. Hence, 
\beq\mathbb B^{p,q}(\Phi,\Phi)= \sum_{c_k>0} c_k \mathbb{B}^{p-k,q-k}(\psi_k,\psi_k)\geq 0.\eeq 
In the following analysis, we will demonstrate that $ m \leq C_{p,q}^k $ holds under appropriate conditions.\\

 \bd
\item If $ q \geq p + 2 $, we have $ s = \min\{p,q-1\} = \min\{p,q\} = t $. By $(1)$ of Lemma \ref{combinatorial2}, 
\beq \min_{0\leq k\leq s} C_{p,q}^k =C_{p,q}^0= \frac{(n-p+1)(p+q)}{2(p+1)}. \eeq
Since $ m \leq \frac{(n-p+1)(p+q)}{2(p+1)}$, we have $ m \leq C_{p,q}^k $ for $ 0 \leq k \leq t$.  Therefore $ \mathbb{B}^{p,q}(\Phi,\Phi)\geq 0 $,  and $\mathbb{B}^{p,q}(\Phi,\Phi)=0$ if and only if all $\psi_k$ are zero. In particular, $\Phi=0$.
\item If $ q = p + 1 $ and $ p \leq \frac{n}{2} $,  we have $ s = t $. By Lemma \ref{combinatorial2},
\beq \min_{0\leq k\leq s} C_{p,q}^k =C_{p,q}^{q-1}= \frac{n+1}{2}. \eeq
Since $ m \leq \frac{n+1}{2}$, we have $ m \leq C_{p,q}^k $ for $ 0 \leq k \leq t $ and therefore $ \mathbb{B}^{p,q} $ is positive.
\item If $ q = p + 1 $ and $ p > \frac{n}{2} $,  $ s = t $. By Lemma \ref{combinatorial2},
\beq \min_{0\leq k\leq s} C_{p,q}^k =C_{p,q}^0= \frac{(n-p+1)(p+q)}{2(p+1)}. \eeq
Since $ m \leq \frac{(n-p+1)(p+q)}{2(p+1)}$, we have $ m \leq C_{p,q}^k $ for $ 0 \leq k \leq t $ and therefore $ \mathbb{B}^{p,q} $ is positive.\ed 

\noindent We shall analyze the semi-positivity of $\mathbb B^{p,q}$.
\bd
\item If $ q \leq p $ and $ p \leq \frac{n}{2} $, we have $ s = q - 1 = t - 1 $. By Lemma \ref{combinatorial2},
\beq \min_{0\leq k\leq s} C_{p,q}^k =C_{p,q}^{q-1}= \frac{n-p+q}{2}. \eeq
Since $ m \leq \frac{n-p+q}{2}$, we have $ m \leq C_{p,q}^k $ for $ 0 \leq k \leq q - 1 $. Moreover, when $ k = q $, by definition one has $ \mathbb{B}^{p-q,0} = 0 $ and so $ \mathbb{B}^{p,q} $ is semi-positive. Suppose that $ \mathbb{B}^{p,q}(\Phi,\Phi) = 0 $, then
\beq \mathbb{B}^{p-k,q-k}(\psi_k,\psi_k) = 0, \eeq
for any $ 0 \leq k \leq t $ and $c_k>0$. For $ 0 \leq k \leq  q-1 $, we have 
\beq 0=\mathbb{B}^{p-k,q-k}(\psi_k,\psi_k) \geq
\frac{4(q-k)(p+1-k)}{p+q-2k}\left(
\sum_{i=1}^m \lambda_i \right)|\psi_k|^2\geq 0.\eeq 
  Hence, we have $ \psi_k = 0 $.  For the last piece $k=q$, we get that  $$ \Phi = L^q \psi $$ for some $ \psi \in \Lambda^{p-q,0} T_x^*M $.
\item If $ q \leq p $ and $ p > \frac{n}{2} $, we have $ s = q - 1 = t - 1 $. By Lemma \ref{combinatorial2},
\beq \min_{0\leq k\leq s} C_{p,q}^k = \frac{(n-p+1)(p+q)}{2(p+1)}. \eeq
Since $ m \leq \frac{(n-p+1)(p+q)}{2(p+1)} $, we have $ m \leq C_{p,q}^k $ for $ 0 \leq k \leq q - 1 $. Moreover, when $ k = q $, one has $ \mathbb{B}^{p-q,0} = 0 $ and therefore $ \mathbb{B}^{p,q} $ is semi-positive. By similar discussions as above, we know that if $ \mathbb{B}^{p,q}(\Phi,\Phi) = 0 $, then $ \Phi = L^q \psi $ for some $ \psi \in \Lambda^{p-q,0} T_x^*M $.
\ed
This completes the proof.
\eproof

\vskip 1\baselineskip 

\noindent \emph{Proof of Theorem \ref{main5}}.
Suppose that $p\leq q$ and $\Phi\in \mathcal H_{\bp}^{p,q}(M,\C)$. By Theorem \ref{main1},
\beq \label{intergralmain2} 0 = (\bar\p_F \Phi, \bar\p_F\Phi) + \frac{1}{4} \int_M \mathbb{B}^{p,q}(\Phi,\Phi) \frac{\omega^n}{n!}. \eeq
\bd \item If $ q \geq p + 2 $ and $m\leq   \frac{(n-p+1)(p+q)}{2(p+1)}$, then by  part $(1)$ of Theorem \ref{positivityofBpq},  $\mathbb{B}^{p,q} $ is positive definite and therefore $\Phi=0$.

\item  Suppose that  $ q = p + 1 $ and $ m\leq \frac{n+1}{2} $. If $p\leq \frac{n}{2}$, then by part $(2)$ of Theorem \ref{positivityofBpq},  $ \mathbb{B}^{p,q} $ is positive definite and  $\Phi=0$. In the case $p>\frac{n}{2}$, by Serre duality, $\mathcal H^{p,q}\cong \mathcal H^{n-q,n-p}$, one can also deduce $\Phi=0$.

\item   Suppose that  $ q = p $ and $ m \leq \frac{n}{2} $. By Serre duality, we can assume that $ p \leq \frac{n}{2}$. By semi-positivity part $(1)$ of Theorem \ref{positivityofBpq}, if $ \mathcal{R} $ is $ \frac{n}{2} $-positive, one has $ \mathbb{B}^{p,p} $ is semi-positive and so $ \mathbb{B}^{p,p}(\Phi,\Phi) = 0 $. Hence, one has $ \Phi= L^p f = f \omega^p $ for some $f\in C^\infty(M)$.  Since $\Phi$ is harmonic, $f$ is constant. \ed
The proof of Theorem \ref{main5} is completed. \qed

\vskip 1\baselineskip

\bproof[Proof of Corollary \ref{cohomologicalspace}] By using Serre duality, we can assume $p\leq \frac{n}{2}$. If $q=p$ or $q=p+1$, then by Theorem \ref{main5}, $\mathcal H^{p,q}(M,\C)=\mathcal H^{p,q}(\C\P^n,\C)$. If $q\geq p+2$, then
\beq  \frac{(n-p+1)(p+q)}{2(p+1)}\geq n-p+1\geq \frac{n}{2}+1>m. \eeq
By part $(1)$ of Theorem \ref{main5}, $\mathcal H^{p,q}(M,\C)=0$.
\eproof

\vskip 1\baselineskip 

\bproof[Proof of Theorem \ref{main7}] Suppose that $ \phi \in
\Omega^{p,q}(M,E) $ is $ \bar\p_E $-harmonic, by Theorem
\ref{KBformulaVBcase}, \beq 0 =
\langle\bar\p{}_F^*\bar\p_F \phi, \phi\rangle +
\frac{1}{4}\mathbb{B}^{p,q}(\phi,\phi) +\left\langle \big(\mathfrak{R}^E \otimes \mathrm{Id}_{\Lambda^{p,q-1}T^*M}\big) \left(\mathbb{S}_\phi\right), \mathbb{S}_\phi \right\rangle. \eeq
Here $ \mathbb{B}^{p,q}: \Omega^{p,q}(M,E) \times \Omega^{p,q}(M,E) \> \mathbb{C} $ is given  by
\beq \mathbb{B}^{p,q}(\psi,\eta) = \left\langle \big(\mathcal{R}\otimes \mathrm{Id}_{\Lambda^{p+1,q-1}T^*M\ts E}\big)\big(\mathbb{T}_{\psi}\big),\mathbb{T}_{\eta} \right\rangle. \eeq
‌It is straightforward to verify that the same conclusion as Theorem \ref{positivityofBpq} holds for the operator $\mathbb B^{p,q}$ defined above.
Intergrating
over $ M $, one has \beq 0 = (\bar\p_F\phi,\bar\p_F\phi) +
\frac{1}{4}\int_M \mathbb{B}^{p,q}(\phi,\phi)\frac {\omega^n}{n!} +
\left( \big(\mathfrak{R}^E \otimes \mathrm{Id}_{\Lambda^{p,q-1}T^*M}\big)
\mathbb{S}_\phi, \mathbb{S}_\phi \right). \eeq Since $ E $ is Nakano
positive, one has $ \mathfrak{R}^E \otimes \mathrm{Id}_{\Lambda^{p,q-1}T^*M}$ is positive-definite. 
\bd
\item If $ p = n $, then $ \mathbb{B}^{p,q} = 0 $ and therefore $ \mathbb{S}_\phi = 0 $. On the other hand,  a straightforward calculation shows $$ |\mathbb{S}_\phi|^2 = q|\phi|^2,$$ and we conclude $ \phi = 0 $. Therefore, $ H_{\bar\p}^{p,q}(M,E) = 0 $ for $q\geq 1$.
\item If $ q \leq p + 1 $, $ p \leq \frac{n}{2} $ and $ m\leq \frac{n-p+q}{2} $. By Theorem \ref{positivityofBpq},  $ \mathbb{B}^{p,q} $ is semi-positive, and one has $ H_{\bar\p}^{p,q}(M,E) = 0 $.
\item If $(p,q)$ is not in the case of $(1)$ or $(2)$ and $ m \leq \frac{(n-p+1)(p+q)}{2(p+1)} $. By Theorem \ref{positivityofBpq}, $ \mathbb{B}^{p,q} $ is also semi-positive, and so $ H_{\bar\p}^{p,q}(M,E) = 0 $.
\ed
This completes the proof.
\eproof

\vskip 2\baselineskip

\section{Weitzenb\"ock formulas with quadratic curvature terms  on  Riemannian manifolds}\label{sectionWR}
In this section we prove Theorem \ref{main8} and establish some applications.  Let $(M,g)$ be a compact and oriented Riemannian manifold. Recall that for any  $\omega\in \Om^{p}(M)$,
$ \mathbb{T}_\omega: \Gamma(M,T^*M)\times  \Gamma(M,T^*M)\>\Om^{p}(M)$ is the operator defined
by: \beq  \mathbb T_\omega(\alpha,\beta) =  \alpha \wedge I_{\beta_\#} \omega - \beta \wedge I_{\alpha_\#} \omega.\eeq
One can also regard it as \beq  \mathbb T_\omega \in \Gamma\left(M,\Lambda^2TM\otimes \Lambda^{p}T^*M\right).\eeq
\btheorem Let $(M,g)$ be a compact Riemannian manifold.
For any   differential form $ \omega \in \Omega^p(M) $, the following Weitzenb\"ock formula  holds
\beq \left\langle\Delta_d \omega, \omega\right\rangle = \left\langle D^*D \omega, \omega\right\rangle + \left\langle \big(\mathfrak{R}\otimes \mathrm{Id}_{\Lambda^pT^*M}\big)\left(\mathbb T_\omega\right),\mathbb T_\omega \right\rangle_{\Lambda^2TM \otimes \Lambda^pT^*M },\label{WeitzenbockR}\eeq
where $D$ is the induced connection on $\Lambda^pT^*M$.
\etheorem

\bproof By using the expression of $d^*$, we have
\beq d^*d \omega
= -g^{jk} I_j \nabla_k(dx^i \wedge \nabla_i\omega).\eeq
A straightforward computation shows
\be
d^*d \omega
& =& -g^{jk} I_j (\nabla_k dx^i \wedge \nabla_i \omega + dx^i \wedge \nabla_k \nabla_i \omega) \\
& =& -g^{jk} (I_j\nabla_k dx^i \wedge \nabla_i \omega - \nabla_k dx^i \wedge I_j \nabla_i \omega) - g^{jk} (\nabla_k \nabla_j \omega - dx^i \wedge I_j \nabla_k \nabla_i \omega) \\
& =& -g^{jk} (\nabla_k \nabla_j \omega - \Gamma_{jk}^i \nabla_i \omega) + g^{jk} (dx^i \wedge I_j \nabla_k \nabla_i \omega - \Gamma_{k\ell}^i dx^\ell \wedge I_j \nabla_i \omega) \\
& =& -g^{jk} (\nabla_k \nabla_j \omega - \Gamma_{jk}^i \nabla_i \omega) + g^{jk} dx^i \wedge I_j(\nabla_k \nabla_i \omega - \Gamma_{ki}^\ell \nabla_\ell \omega).  \ee
On the other hand, if we write $s=\omega\in\Gamma(M,E)$ where $E=\Lambda^pT^*M$,
\beq D^*D s
= -g^{jk} I_j \nabla^{T^*M \otimes E}_k(dx^i \wedge  \nabla^E_is). \eeq
Therefore,
\be
D^*D s
& =& -g^{jk} I_j (\nabla_k dx^i \wedge \nabla^E_i s + dx^i \wedge \nabla^E_k \nabla^E_i s) \\
& =& -g^{jk} (-\Gamma^i_{jk} \nabla^E_is + \nabla^E_k \nabla^E_j s) \\
& =& -g^{jk} (\nabla_k \nabla_j \omega - \Gamma^i_{jk} \nabla_i\omega).
\ee

\noindent  We define $ \alpha \in \Gamma(M,T^*M \otimes T^*M \otimes \wedge^{p-1} T^*M) $ as
\beq \alpha(X,Y) = I_X \nabla_Y \omega.  \eeq
It is easy to see that \be
(\nabla_Z \alpha)(X,Y)
& =& \nabla_Z (\alpha(X,Y)) - \alpha(\nabla_Z X, Y) - \alpha(X, \nabla_Z Y) \\
& =& \nabla_Z I_X \nabla_Y \omega - I_{\nabla_Z X} \nabla_Y \omega - I_X \nabla_{\nabla_Z Y} \omega \\
& = & I_X \nabla_Z \nabla_Y \omega - I_X \nabla_{\nabla_Z Y} \omega
\ee
where the third identity follows from the contraction formula
\beq \nabla_Z\circ
I_W=I_W\circ \nabla_ZZ+ I_{\nabla_ZW}. \eeq
 Therefore, we obtain
\be
dd^* \omega
=- dx^i \wedge \nabla_i(g^{jk}I_j \nabla_k \omega)
&=&- dx^i \wedge \nabla_i (\mathrm{tr}_g \alpha)\\
& =& - dx^i \wedge g^{jk} (\nabla_i  \alpha)\left(\frac{\p}{\p x^j},\frac{\p}{\p x^k}\right) \\
& = &- dx^i \wedge g^{jk} I_j(\nabla_i \nabla_k \omega - \Gamma_{ki}^\ell \nabla_\ell \omega).
\ee
Now we obtain
\beq
\Delta_d \omega - D^*D \omega
= g^{jk} dx^i \wedge I_j(\nabla_k \nabla_i \omega - \nabla_i \nabla_k \omega)= g^{jk}dx^i \wedge I_k \left(R\left( \frac{\p}{\p x^j}, \frac{\p}{\p x^i} \right) \omega\right).
\eeq
On the other hand, at a fixed point with orthonormal frame, it is easy to deduce that
\beq \mathbb T_\omega=\sum_{i<j} \left(\frac{\p}{\p x^i}\wedge \frac{\p}{\p x^j}\right)\ts \left(dx^i \wedge I_j \omega - dx^j \wedge I_i \omega\right),\eeq
and
\be &&\left\langle \big(\mathfrak{R}\otimes \mathrm{Id}_{\Lambda^pT^*M}\big)\left(\mathbb T_\omega\right),\mathbb T_\omega \right\rangle_{\Lambda^2TM \otimes \Lambda^pT^*M }\\&= & \sum_{i < j} \sum_{k < \ell} R_{ij\ell k} \langle dx^i \wedge I_j \omega - dx^j \wedge I_i \omega, dx^k \wedge I_\ell \omega - dx^\ell \wedge I_k \omega \rangle \\
&= & \sum_{i, j, k, \ell } R_{ij\ell k} \langle dx^i \wedge I_j \omega , dx^k \wedge I_\ell \omega \rangle \\
&= & -\sum_{i, j, \ell} \left\langle dx^i \wedge I_j \omega ,\left( R \left( \frac{\p}{\p x^j},\frac{\p}{\p x^i} \right) dx^\ell \right) \wedge I_\ell \omega \right\rangle  \\
&=& -\sum_{i, j, \ell} \left\langle dx^i \wedge I_j \omega ,R \left( \frac{\p}{\p x^j},\frac{\p}{\p x^i} \right) \omega \right\rangle.
\ee
This is exactly $\left\langle \Delta_d \omega - D^*D \omega,\omega\right\rangle$ and we complete the proof of \eqref{WeitzenbockR}.
\eproof

\btheorem \label{PWE}
For any $ \omega\in \Omega^k(M) $,
\beq  |\mathbb T_\omega|^2 \leq 2\min\{k,n-k\} |\omega|^2. \eeq
\etheorem
\noindent  Theorem \ref{PWE} is established in \cite{PetersenWink2021a}. The proof of it follows by a formal linear algebra using the operator $\mathbb T_\omega$ which is  simpler than that of Theorem \ref{main3}.  Let $V$ be a real vector
space and $\dim_\R V=r$. Let $A:V\>V$ be a linear map. For any $1\leq
p\leq r$, the $p$-th compound matrix  $\wedge^pA\in\mathrm{End}(\Lambda^pV)$ of $A$ is defined as \beq
\left(\wedge^pA\right)(v_1\wedge \cdots \wedge v_p)=\sum_{i=1}^p
v_1\wedge \cdots\wedge A v_i\wedge \cdots\wedge v_p.\eeq
Since  $ \mathbb T_\omega \in \Gamma\left(M,\Lambda^2TM\otimes \Lambda^{p}T^*M\right)$,
for a given vector $v\in \Lambda^2TM$,  we define a map $T( v, \bullet): \Om^p(M)\>\Om^p(M)$ by
\beq T(v,\omega): =\mathbb T_\omega (v)\in \Om^p(M),\eeq
and this map is denoted by $T_p(v)$. One can see clearly that
\beq  \wedge^k T_1(v)=T_k(v).\eeq
At a fixed point $x\in M$, we assume that $\{e_i=\frac{\p}{\p x^i}\}$ is an orthonormal frame of $T_xM$.  Let $v=\frac{1}{2}\sum v^{ij} e_i\wedge e_j\in\wedge^2 T_xM$ where $v^{ij}$ is skew-symmetric. Suppose the matrix $V:=\left[ v^{ij}\right]$ has complex eigenvalues $\lambda_1,\cdots, \lambda_n$. Since $V$ is skew-symmetric, its rank is denoted by $2r$. We also assume that \beq (\lambda_1,\cdots, \lambda_n)=\left(\sq \Lambda_1,\cdots, \sq \Lambda_r, -\sq \Lambda_1,\cdots, -\sq \Lambda_r,0,\cdots,0\right)\eeq
where  $\Lambda_1\geq \cdots \geq \Lambda_r>0$. In this case, $\lambda_{2r+1}=\cdots=\lambda_n=0$. For simplicity, we also set $\Lambda_{r+1}=\cdots=\Lambda_{n}=0$. Hence, one obtains:

\blemma\label{compoundeigenvalue}\bd \item  The eigenvalues of  $T_1(v):T^*_xM\>T_x^*M$ are $2\lambda_1,\cdots, 2\lambda_n$.
\item The eigenvalues of  $T_k(v)$ are of the form $2\sum\limits_{s=1}^k \lambda_{i_s}$ where $1\leq i_1<\cdots<i_k\leq n$.
\ed  In particular, the maximal absolute value  of the  eigenvalues of  $T_k(v)$ is \beq \min\{2\left(\Lambda_1+\cdots +\Lambda_k\right), 2\left(\Lambda_1+\cdots +\Lambda_{n-k}\right) \}.\eeq
\elemma
\noindent
Indeed, let $\lambda$ be an eigenvalue of $T_1(v)$, $a=\sum a_s dx^s$ be an eigenvector and  $$T_1(v)(a)=\lambda a.$$ A straightforward computation shows
\beq T_1(v)(a)=\sum_{i,j} v^{ij}a_s\left( g_{ im} dx^m I_j(dx^s) - g_{jm} dx^mI_i(dx^s) \right)=2\sum_{j,m} v^{mj}a_j dx^m.\eeq
In particular,  one has
\beq 2\sum_{j} v^{mj}a_j=\lambda a_m.\eeq
Hence,  $\lambda=2\lambda_i$ for some $i$.

\bproof[Proof of Theorem \ref{PWE}]  Let $v=\frac{1}{2}\sum v^{ij} e_i\wedge e_j\in\wedge^2 T_xM$ where $v^{ij}$ is skew-symmetric.
It is easy to see that
\beq |v|^2_g=\sum_{i,j} |v^{ij}|^2=-\sum_{i,j}v^{ij}v^{ji}=- \mathrm{tr}\left(V\circ V\right)=-\sum_i \lambda_i^2=2\sum_{i}\Lambda_i^2.\eeq
Hence, by Lemma \ref{compoundeigenvalue},
\be   |\mathbb T_\omega(v)|^2=|T_k(v)(w)|^2&\leq& \min\{4\left(\Lambda_1+\cdots +\Lambda_k\right)^2, 4\left(\Lambda_1+\cdots +\Lambda_{n-k}\right)^2 \} |w|^2\\ &\leq&2\min\{k,n-k\}  |\omega|^2 |v|_g^2,\ee
where the last step follows from  Cauchy-Schwarz inequalities.
\eproof

\noindent As an application of Theorem \ref{main8} and Theorem \ref{PWE}, one obtains
the following result of  Petersen-Wink \cite{PetersenWink2021a}:
\btheorem Let $(M^n,g)$ be a compact   Riemannian manifold. \bd

\item If $ \mathfrak{R} $ is $ p $-positive, then $ b_1(M) = \cdots = b_{n-p}(M) = 0$ and $b_p(M) = \cdots= b_{n-1}(M) = 0 $.

\item If $\mathfrak{R} $ is $ p $-positive and $ 2p \leq n $, then $ b_k(M) = 0 $ for all $ 1 \leq k \leq n - 1$.

\item If  $ \mathfrak{R} $ is $p$-semipositive, then any harmonic $k$-form is parallel for $ 1 \leq k \leq n - p $, or $ p \leq k \leq n-1 $.

\ed
\etheorem

\vskip 2\baselineskip

\section{Weitzenb\"ock formulas with quadratic curvature terms on K\"ahler manifolds}\label{sectionWK}

The Weitzenb\"ock formulas with quadratic curvature terms on compact K\"ahler manifolds constitute a natural complexification of their Riemannian counterparts. Here we present the principal results; their proofs follow straightforwardly from the established framework. Let $(M,\omega_g) $ be a compact K\"ahler manifold. The reduced (complexified)  curvature operator $ \mathscr{R}: \Gamma(M,T^{1,0}M \otimes T^{0,1}M) \> \Gamma(M,T^{1,0}M \otimes T^{0,1}M)$ is
defined as:
\beq \left\langle \mathscr{R}\left(\frac{\p}{\p z^i} \wedge \frac{\p}{\p \bar z{}^j}\right), \frac{\p}{\p z^\ell} \wedge \frac{\p}{\p \bar z{}^k}\right\rangle = R_{i\bar j k\bar\ell}. \eeq
For any $ \phi \in \Omega^{p,q}(M) $, $ \mathbb{Y}_\phi: \Gamma(M,T^{*1,0}M) \times \Gamma(M,T^{*0,1}M) \> \Omega^{p,q}(M) $ is defined
by:
\beq \mathbb{Y}_\phi(\alpha,\beta) = \beta \wedge I_{\alpha_\#} \phi - \alpha \wedge I_{\beta_\#} \phi. \eeq
 It is obvious that
\beq  \mathbb{Y}_\phi \in \Gamma\left(M,\left(T^{1,0}M \wedge T^{0,1}M\right) \otimes \Lambda^{p,q}T^*M\right). \eeq

\btheorem Let $(M,\omega_g)$ be a compact K\"ahler manifold. For any differential form $ \phi \in \Omega^{p,q}(M) $, the following Weitzenb\"ock formula  holds
\beq \langle\Delta_d \phi, \phi\rangle = \langle D^*D\phi, \phi\rangle + \left\langle \big(\mathscr{R}\otimes \mathrm{Id}_{\Lambda^{p,q}T^*M}\big)\big(\mathbb Y_\phi\big),\mathbb Y_\phi \right\rangle_{\left(T^{1,0}M \wedge T^{0,1}M\right) \otimes \Lambda^{p,q}T^*M }, \eeq
where $ D $ is the induced connection on $ \Lambda^{p,q}T^*M $.
\etheorem

\noindent The following results are essentially established in  \cite{PetersenWink2021b}.
\btheorem
Given  $ \phi \in \Omega^{p,q}(M) $, if there exists $ k \geq 0 $ and $ \psi \in \Omega^{(p-k),(q-k)}(M)$ such that $ \phi = L^k \psi $ and $ \Lambda \psi = 0 $, the following inequality holds
\beq |\mathbb{Y}_\phi|^2 \leq (p+q-2k)|\phi|^2. \eeq

\etheorem

\btheorem Let $(M,\omega_g)$ be a compact K\"ahler manifold. Assume that the reduced curvature operator $ \mathscr{R}: \Gamma(M,T^{1,0}M \wedge T^{0,1}M) \> \Gamma(M,T^{1,0}M \wedge T^{0,1}M) $ is $m$-positive.
\bd
\item If $ p \neq q $ and $ m \leq \left(n+1 - \frac{p^2+q^2}{p+q}\right) $, one has
\beq H^{p,q}_{\bar \p}(M,\mathbb{C}) = 0. \eeq
\item If $ m \leq n+1-p $, one has
\beq H^{p,p}_{\bar \p}(M,\mathbb{C}) = \mathbb{C}. \eeq
\ed
\etheorem

\vskip 2\baselineskip


\begin{thebibliography}{BCHM}




    \bibitem[BG24]{BettiolGoodman24}
    Renato Bettiol and McFeely Jackson Goodman.
    \newblock Curvature operators and rational cobordism.
    \newblock {\em Adv. Math.}, 458,2024.

    \bibitem[Boc46]{Bochner46}
    Salomon Bochner.  Vector fields and Ricci curvature, {\emph{Bull. Amer. Math. Soc}}. 52: 776--797,  1946
    \bibitem[Boc48]{Bochner48}
    Salomon Bochner.
    \newblock Curvature and Betti numbers.
    \newblock {\em Ann. of Math.}, 49(2):379--390, 1948.


    \bibitem[Boc49]{Bochner49}
    Salomon Bochner.   Curvature and Betti numbers II.  \newblock {\em Ann. of Math.}, 50(2):77--93, 1949.

    \bibitem[BW08]{BohmWilking2008}
    Christoph B\"ohm and Burkhard Wilking.
    \newblock Manifolds with positive curvature operators are space forms.
    \newblock {\em Ann. of Math.}, 167(2):1079--1097, 2008.



    \bibitem[Bre10]{Brendle2010}
    Simon Brendle.
    \newblock Einstein manifolds with nonnegative isotropic curvature are locally symmetric.
    \newblock {\em Duke Math. J.}, 151(1):1--21, 2010.

    \bibitem[Bre19]{Brendle2019}
    Simon Brendle.
    \newblock Ricci flow with surgery on manifolds with positive isotropic curvature.
    \newblock {\em Ann. of Math.}, 190(2):465--559, 2019.

    \bibitem[BS08]{BrendleSchoen2008}
    Simon Brendle and Richard Schoen.
    \newblock Classification of manifolds with weakly $1/4$-pinched curvatures.
    \newblock {\em Acta Math.},200:1--13,2008;

    \bibitem[BS09]{BrendleSchoen2009}
    Simon Brendle and Richard Schoen.
    \newblock Manifolds with 1/4-pinched curvature are space forms.
    \newblock {\em J. Am. Math. Soc.}, 22(1):287--307, 2009.


\bibitem[BNPSW25]{BNPSW25}  Kyle Broder, Jan Nienhaus, Peter Petersen, James Stanfield and Matthias Wink. Vanishing theorems for Hodge numbers and the Calabi curvature operator.   arXiv:2503.06870


    \bibitem[CD23]{ChoDung2023}
    Gunhee Cho and Nguyen Thac Dung.
    \newblock Vanishing results from Lichnerowicz Laplacian on complete K\"ahler manifolds and applications.
    \newblock {\em J. Math. Anal. Appl.},517,2023.



    \bibitem[CTZ12]{ChenTangZhu2012}
    Bing-Long Chen, Siu-Hung Tang and Xi-Ping Zhu.
    \newblock Complete classification of compact four-manifolds with positive isotropic curvature.
    \newblock {\em J. Differ. Geom.}, 91(1):41--80, 2012.

    \bibitem[CZ06]{ChenZhu2006}
    Bing-Long Chen and Xi-Ping Zhu.
    \newblock Ricci flow with surgery on four-manifolds with positive isotropic curvature.
    \newblock {\em J. Differ. Geom.}, 74(2):177--264, 2006.







    \bibitem[CGT23]{CGT23}Xiaodong Cao, Matthew J.  Gursky and Hung Tran.
    Curvature of the second kind and a conjecture of Nishikawa. {\emph{Comment. Math. Helv.}},  98: 195--216, 2023.


\bibitem[CV60]{CV60}] Eugenio Calabi and Edoardo Vesentini.  On compact, locally symmetric K\"ahler manifolds,  \emph{Ann.
of Math. (2)}  71 , 472--507, 1960.


    \bibitem[DF24]{DaiFu2024}
    Zhi-Lin Dai and Hai-Ping Fu.
    \newblock Einstein manifolds and curvature operator of the second kind.
    \newblock {\em Calc. Var. Partial Differential Equations}, 63, 2024.

    \bibitem[ES64]{EellsSampson1964}
    Jams Eells and Joseph H. Sampson. Harmonic mappings of Riemannian manifolds. \emph{Amer. J. Math.}  86:109--160, 1964.



    \bibitem[GM75]{GM75}
    Sylvestre
    Gallot and Daniel Meyer.
    Opérateur de courbure et laplacien des formes différentielles d'une variété riemannienne. {\em{J. Math. Pures Appl.}}, 54: 259--284, 1975








    \bibitem[Ham82]{Hamilton1982}
    Richard S. Hamilton.
    \newblock Three-manifolds with positive Ricci curvature.
    \newblock {\em J. Differ. Geom.}, 17:255--306, 1982.

    \bibitem[Ham86]{Hamilton1986}
    Richard S. Hamilton.
    \newblock Four-manifolds with positive curvature operator.
    \newblock {\em J. Differ. Geom.}, 24:153--179, 1986.

\bibitem[Huy05]{Huy05} Daniel Huybrechts, {\it Complex geometry}, Universitext, Springer, Berlin, 2005.

\bibitem[HJ85]{HJ85}
    Roger~A. Horn and Charles~R. Johnson. {\it Matrix analysis}, Cambridge Univ. Press, Cambridge, 1985.

    \bibitem[Li23]{Li23} Xiaolong Li.
    K\"ahler manifolds and the curvature operator of the second kind. {\emph{Math. Z.}},  303: Paper No. 101, 26 pp, 2023.

    \bibitem[Li24]{Li24} Xiaolong Li. Manifolds with nonnegative curvature operator of the second kind, \emph{ Commun. Contemp. Math.}  26: Paper No. 2350003, 26 pp, 2024.

     \bibitem[LY12]{LY12} Kefeng Liu and Xiaokui Yang. Geometry of Hermitian manifolds.
    \emph{Internat. J. Math.} \textbf{23}, 40pp., 2012.



    \bibitem[Mey71]{Meyer1971}
    Daniel  Meyer.
    \newblock Sur les vari\'{e}t\'{e}s riemanniennes \`{a} op\'{e}rateur de courbure positif.
    \newblock {\em Comptes Rendus Acad. Sci. Paris S\'{e}r. A--B}, 272 A482--A485, 1971.

    \bibitem[MM88]{MicallefMoore1988}
    Mario J. Micallef and John D. Moore.
    \newblock Minimal two-spheres and the topology of manifolds with positive curvature on totally isotropic two-planes.
    \newblock {\em Ann. of Math.}, 127(1):199--227, 1988.



    \bibitem[Mori79]{Mori1979}
    Shigefumi Mori.
    \newblock Projective manifolds with ample tangent bundles.
    \newblock {\em Ann. of Math.}, 110(3):593--606, 1979.

    \bibitem[NPW23]{NienhausPetersenWink2023}
    Jan  Nienhaus, Peter Petersen and Matthias Wink.
    \newblock Betti numbers and the curvature operator of the second kind.
    \newblock {\em J. Lond. Math. Soc.},108(2):1642--1668, 2023.





    \bibitem[Pet16]{Pet16} Peter Petersen. {\it Riemannian geometry}, $3$rd edition,
    Graduate Texts in Mathematics, 171, 2016.

    \bibitem[PW21a]{PetersenWink2021a}
    Peter Petersen and Matthias Wink.
    \newblock New curvature conditions for the Bochner technique.
    \newblock {\em Invent. math.}, 224:33--54, 2021.

    \bibitem[PW21b]{PetersenWink2021b}
    Peter Petersen and Matthias Wink.
    \newblock Vanishing and estimation results for Hodge numbers.
    \newblock {\em J. reine angew. Math.}, 780:197--219, 2021.

    \bibitem[PW22]{PetersenWink2022}
    Peter Petersen and Matthias Wink.
    \newblock Tachibana-type theorems and special holonomy.
    \newblock {\em Ann. Global Anal. Geom.},61:847--868,2022.


    \bibitem[SY80]{SiuYau1980}
    Yum-Tong Siu and Shing-Tung Yau.
    \newblock Compact K\"ahler manifolds of positive bisectional curvature.
    \newblock {\em Invent. math.}, 59(2):189--204, 1980.

    \bibitem[Siu80]{Siu1980}
    Yum-Tong Siu.   The complex-analyticity of harmonic maps and the strong rigidity of compact K\"ahler manifolds,  {\em Ann. of Math.}  112(2): 73--111, 1980.


    \bibitem[Siu82]{Siu1982}
    Yum-Tong Siu.
    Complex-analyticity of harmonic maps, vanishing and Lefschetz theorems.  \emph{J. Differential Geom.} 17:  55--138, 1982.


    \bibitem[Wil13]{Wil13} Burkhard Wilking.
    A Lie algebraic approach to Ricci flow invariant curvature conditions and Harnack inequalities. {\emph{J. Reine Angew. Math.}}, 679: 223--247, 2013.


    \bibitem[Wu88] {Wu88}
    Hung-Hsi Wu. {\it The Bochner technique in differential geometry.}  Math. Rep. {3}(2): 289--538, 1988. 
    
    
    \bibitem[Xu25]{Xu2025}
    Kai  Xu.
    \newblock Dimension constraints in some problems involving intermediate curvature.
    \newblock {\em Trans. Amer. Math. Soc.}, 378:2091--2112,2025.

    \bibitem[Yang18]{Yang2018}
    Xiaokui Yang.
    \newblock RC-positivity, rational connectedness and Yau's conjecture.
    \newblock {\em Camb. J. Math.}, 6(2):183--212, 2018.

  \bibitem[YZ25]{YZ25}
 Xiaokui Yang and Liangdi Zhang.
New curvature characterizations for spherical space forms and complex projective spaces.
\emph{Trans. Amer. Math. Soc.} 378 (1): 679--694, 2025.



\end{thebibliography}

\end{document}